\documentclass[reqno,11pt]{amsart}
\usepackage{amsmath}
\usepackage{amsxtra}
\usepackage{amscd}
\usepackage{amsthm}
\usepackage{amsfonts}
\usepackage{amssymb}
\usepackage{eucal}
\usepackage[matrix,arrow,curve]{xy}

\theoremstyle{plain}
\newtheorem{Thm}[subsection]{Theorem}
\newtheorem{Cor}[subsection]{Corollary}
\newtheorem{Lem}[subsection]{Lemma}
\newtheorem{Prop}[subsection]{Proposition}
\newtheorem{Conj}[subsection]{Conjecture}

\theoremstyle{definition}
\newtheorem{Def}[subsection]{Definition}

\theoremstyle{remark}

\newtheorem{Rem}[subsection]{Remark}

\newtheorem{conj}{Conjecture}

\numberwithin{equation}{section}
\renewcommand{\rm}{\normalshape}

\newif\ifShowLabels
\ShowLabelstrue
\newdimen\theight
\def\TeXref#1{%
    \leavevmode\vadjust{\setbox0=\hbox{{\tt
        \quad\quad  {\small \rm #1}}}%
    \theight=\ht0
    \advance\theight by \lineskip
    \kern -\theight \vbox to
    \theight{\rightline{\rlap{\box0}}%
    \vss}%
    }}%

\ShowLabelsfalse

\renewcommand{\sec}[2]{\section{#2}\label{S:#1}%
    \ifShowLabels \TeXref{{S:#1}} \fi}
\newcommand{\ssec}[2]{\subsection{#2}\label{SS:#1}%
    \ifShowLabels \TeXref{{SS:#1}} \fi}

\newcommand{\refs}[1]{Section ~\ref{S:#1}}
\newcommand{\refss}[1]{Section ~\ref{SS:#1}}

\newcommand{\reft}[1]{Theorem ~\ref{T:#1}}
\newcommand{\refl}[1]{Lemma ~\ref{L:#1}}
\newcommand{\refp}[1]{Proposition ~\ref{P:#1}}
\newcommand{\refc}[1]{Corollary ~\ref{C:#1}}

\newcommand{\refe}[1]{\eqref{E:#1}}

\newenvironment{thm}[1]%
    { \begin{Thm} \label{T:#1}  \ifShowLabels \TeXref{T:#1} \fi }%
    { \end{Thm} }

\renewcommand{\th}[1]{\begin{thm}{#1} \sl }
\renewcommand{\eth}{\end{thm} }

\newenvironment{lemma}[1]%
    { \begin{Lem} \label{L:#1}  \ifShowLabels \TeXref{L:#1} \fi }%
    { \end{Lem} }
\newcommand{\lem}[1]{\begin{lemma}{#1} \sl}
\newcommand{\elem}{\end{lemma}}

\newenvironment{propos}[1]%
    { \begin{Prop} \label{P:#1}  \ifShowLabels \TeXref{P:#1} \fi }%
    { \end{Prop} }
\newcommand{\prop}[1]{\begin{propos}{#1}\sl }
\newcommand{\eprop}{\end{propos}}

\newenvironment{corol}[1]%
    { \begin{Cor} \label{C:#1}  \ifShowLabels \TeXref{C:#1} \fi }%
    { \end{Cor} }
\newcommand{\cor}[1]{\begin{corol}{#1} \sl }
\newcommand{\ecor}{\end{corol}}

\newenvironment{defeni}[1]%
    { \begin{Def} \label{D:#1}  \ifShowLabels \TeXref{D:#1} \fi }%
    { \end{Def} }
\newcommand{\defe}[1]{\begin{defeni}{#1} \sl }
\newcommand{\edefe}{\end{defeni}}

\newenvironment{remark}[1]%
    { \begin{Rem} \label{R:#1}  \ifShowLabels \TeXref{R:#1} \fi }%
    { \end{Rem} }
\newcommand{\rem}[1]{\begin{remark}{#1}}
\newcommand{\erem}{\end{remark}}

\newenvironment{conjec}[1]%
    { \begin{Conj} \label{Co:#1}  \ifShowLabels \TeXref{Co:#1} \fi }%
    { \end{Conj} }
\renewcommand{\conj}[1]{\begin{conjec}{#1} \sl }
\newcommand{\econj}{\end{conjec}}

\newcommand{\eq}[1]%
    { \ifShowLabels \TeXref{E:#1} \fi
       \begin{equation} \label{E:#1} }
\newcommand{\eeq}{ \end{equation} }

\newcommand{\prf}{ \begin{proof} }
\newcommand{\epr}{ \end{proof} }

\newcommand\alp{\alpha}     

\newcommand\gam{\gamma}     \newcommand\Gam{\Gamma}

\newcommand\kap{\kappa}
\newcommand\lam{\lambda}        \newcommand\Lam{\Lambda}

\newcommand\calB{{\mathcal{B}}}

\newcommand\calD{{\mathcal{D}}}

\newcommand\calF{{\mathcal{F}}}

\newcommand\calJ{{\mathcal{J}}}
\newcommand\calK{{\mathcal{K}}}
\newcommand\calL{{\mathcal{L}}}
\newcommand\calM{{\mathcal{M}}}

\newcommand\calO{{\mathcal{O}}}

\newcommand\calU{{\mathcal{U}}}
\newcommand\calV{{\mathcal{V}}}
\newcommand\calW{{\mathcal{W}}}
\newcommand\calX{{\mathcal{X}}}

        \newcommand\bfC{{\mathbf C}}

        \newcommand\bfU{{\mathbf U}}

\newcommand\PP{\mathbb{P}}
\renewcommand\AA{\mathbb{A}}

\newcommand\GG{\mathbb{G}}

\newcommand\ZZ{\mathbb{Z}}

\newcommand\CC{\mathbb{C}}

 \newcommand\grg{{\mathfrak{g}}}

\newcommand\sdp{\times \hskip -0.3em {\raise 0.3ex
\hbox{$\scriptscriptstyle |$}}} 

\newcommand\Gr{\operatorname{Gr}}

\newcommand\Int{\operatorname{Int}}

\newcommand\SL{{\rm SL}}

\newcommand\Spec{\operatorname{Spec}}

\newcommand\Sym{\operatorname{Sym}}

\newcommand\of{{\overline{f}}}

\newcommand\tilZ{{\widetilde{Z}}}

\newcommand\x{\times}
\newcommand\ten{\otimes}

\newcommand{\ra}{\rangle}
\newcommand{\la}{\langle}

\renewcommand\Spec{\operatorname{Spec}}

\newcommand\nc{\newcommand}

\newcommand{\go}{G(\calO)}

\newcommand{\ogg}{{\overline \gg}}
\newcommand{\ogl}{\ogg^{\lam}}

\newcommand{\iso}{{\stackrel{\sim}{\longrightarrow}}}

\nc\aff{\operatorname{aff}}
\nc\oGr{\overline{\Gr}}
\nc\Bun{\operatorname{Bun}}
\nc\hgrg{\widehat{\grg}}
\renewcommand\Int{\operatorname{Int}}
\nc\bInt{\overline{\Int}}
\nc\hatLam{\widehat{\Lam}}
\nc\bmu{\overline{\mu}}
\nc\bnu{\overline{\nu}}
\nc\blambda{\overline{\lam}}
\renewcommand\SL{\operatorname{SL}}
\nc\ocalW{\overline{\calW}}
\nc\pos{\operatorname{pos}}
\nc\IH{\operatorname{IH}}
\nc\Rep{\operatorname{Rep}}
\nc\Gal{\operatorname{Gal}}
\nc{\tilGr}{\widetilde{\Gr}}

\nc\Pic{\operatorname{Pic}}

\nc{\HC}{{\mathcal{HC}}}
\nc{\on}{\operatorname}
\nc{\BA}{{\mathbb{A}}}
\nc{\BC}{{\mathbb{C}}}
\nc{\BG}{{\mathbb{G}}}
\nc{\BM}{{\mathbb{M}}}
\nc{\BN}{{\mathbb{N}}}
\nc{\BQ}{{\mathbb{Q}}}
\nc{\BP}{{\mathbb{P}}}
\nc{\BR}{{\mathbb{R}}}
\nc{\BZ}{{\mathbb{Z}}}
\nc{\BS}{{\mathbb{S}}}

\nc{\CA}{{\mathcal{A}}}
\nc{\CB}{{\mathcal{B}}}
\nc{\CalC}{{\mathcal C}}
\nc{\CalD}{{\mathcal D}}
\nc{\CE}{{\mathcal{E}}}
\nc{\CF}{{\mathcal{F}}}
\nc{\CG}{{\mathcal{G}}}
\nc{\CH}{{\mathcal{H}}}
\nc{\CK}{{\mathcal{K}}}
\nc{\CL}{{\mathcal{L}}}
\nc{\CM}{{\mathcal{M}}}
\nc{\CMM}{{\mathcal{M}^{\operatorname{gen}}_\hbar(-\rho)}}
\nc{\CN}{{\mathcal{N}}}
\nc{\CO}{{\mathcal{O}}}
\nc{\CP}{{\mathcal{P}}}
\nc{\CQ}{{\mathcal{Q}}}
\nc{\CR}{{\mathcal{R}}}
\nc{\CS}{{\mathcal{S}}}
\nc{\CT}{{\mathcal{T}}}
\nc{\CU}{{\mathcal{U}}}
\nc{\CV}{{\mathcal{V}}}
\nc{\CW}{{\mathcal{W}}}
\nc{\CX}{{\mathcal{X}}}
\nc{\CY}{{\mathcal{Y}}}
\nc{\CZ}{{\mathcal{Z}}}

\nc{\gen}{{\operatorname{gen}}}
\nc{\cM}{{\check{\mathcal M}}{}}
\nc{\csM}{{\check{\mathcal A}}{}}
\nc{\obM}{{\overset{\circ}{\mathbf M}}{}}
\nc{\oCA}{{\overset{\circ}{\mathcal A}}{}}
\nc{\obA}{{\overset{\circ}{\mathbf A}}{}}
\nc{\ooM}{{\overset{\circ}{M}}{}}
\nc{\osM}{{\overset{\circ}{\mathsf M}}{}}
\nc{\vM}{{\overset{\bullet}{\mathcal M}}{}}
\nc{\nM}{{\underset{\bullet}{\mathcal M}}{}}
\nc{\obD}{{\overset{\circ}{\mathbf D}}{}}
\nc{\cp}{{\overset{\circ}{\mathbf p}}{}}
\nc{\ofZ}{{\overset{\circ}{\mathfrak Z}}{}}

\nc{\fa}{{\mathfrak{a}}}
\nc{\fb}{{\mathfrak{b}}}
\nc{\fg}{{\mathfrak{g}}}
\nc{\fgl}{{\mathfrak{gl}}}
\nc{\fh}{{\mathfrak{h}}}
\nc{\fj}{{\mathfrak{j}}}
\nc{\fm}{{\mathfrak{m}}}
\nc{\fn}{{\mathfrak{n}}}
\nc{\fu}{{\mathfrak{u}}}
\nc{\fp}{{\mathfrak{p}}}
\nc{\frr}{{\mathfrak{r}}}
\nc{\fs}{{\mathfrak{s}}}
\nc{\ft}{{\mathfrak{t}}}
\nc{\fT}{{\mathfrak{T}}}
\nc{\ofT}{{\overline{\mathfrak T}}}
\nc{\ofS}{{\overline{\mathfrak S}}}
\nc{\fsl}{{\mathfrak{sl}}}
\nc{\hsl}{{\widehat{\mathfrak{sl}}}}
\nc{\hgl}{{\widehat{\mathfrak{gl}}}}
\nc{\hg}{{\widehat{\mathfrak{g}}}}
\nc{\chg}{{\widehat{\mathfrak{g}}}{}^\vee}
\nc{\hn}{{\widehat{\mathfrak{n}}}}
\nc{\chn}{{\widehat{\mathfrak{n}}}{}^\vee}

\nc{\fA}{{\mathfrak{A}}}
\nc{\fB}{{\mathfrak{B}}}
\nc{\fD}{{\mathfrak{D}}}
\nc{\fE}{{\mathfrak{E}}}
\nc{\fF}{{\mathfrak{F}}}
\nc{\fG}{{\mathfrak{G}}}
\nc{\fI}{{\mathfrak{I}}}
\nc{\fJ}{{\mathfrak{J}}}
\nc{\fK}{{\mathfrak{K}}}
\nc{\fL}{{\mathfrak{L}}}
\nc{\fM}{{\mathfrak{M}}}
\nc{\fN}{{\mathfrak{N}}}
\nc{\frP}{{\mathfrak{P}}}
\nc{\fS}{{\mathfrak S}}
\nc{\fU}{{\mathfrak{U}}}
\nc{\fZ}{{\mathfrak{Z}}}

\nc{\bb}{{\mathbf{b}}}
\nc{\bc}{{\mathbf{c}}}
\nc{\be}{{\mathbf{e}}}
\nc{\bj}{{\mathbf{j}}}
\nc{\bn}{{\mathbf{n}}}
\nc{\bp}{{\mathbf{p}}}
\nc{\bq}{{\mathbf{q}}}
\nc{\bv}{{\mathbf{v}}}
\nc{\bx}{{\mathbf{x}}}
\nc{\by}{{\mathbf{y}}}
\nc{\bw}{{\mathbf{w}}}
\nc{\bA}{{\mathbf{A}}}
\nc{\bB}{{\mathbf{B}}}
\nc{\bC}{{\mathbf{C}}}
\nc{\bK}{{\mathbf{K}}}
\nc{\bD}{{\mathbf{D}}}
\nc{\bH}{{\mathbf{H}}}
\nc{\bM}{{\mathbf{M}}}
\nc{\bN}{{\mathbf{N}}}
\nc{\bS}{{\mathbf{S}}}
\nc{\bT}{{\mathbf{T}}}
\nc{\bV}{{\mathbf{V}}}
\nc{\bW}{{\mathbf{W}}}
\nc{\bX}{{\mathbf{X}}}
\nc{\bP}{{\mathbf{P}}}
\nc{\bZ}{{\mathbf{Z}}}

\nc{\sA}{{\mathsf{A}}}
\nc{\sB}{{\mathsf{B}}}
\nc{\sC}{{\mathsf{C}}}
\nc{\sD}{{\mathsf{D}}}
\nc{\sF}{{\mathsf{F}}}
\nc{\sK}{{\mathsf{K}}}
\nc{\sM}{{\mathsf{M}}}
\nc{\sO}{{\mathsf{O}}}
\nc{\sQ}{{\mathsf{Q}}}
\nc{\sP}{{\mathsf{P}}}
\nc{\sV}{{\mathsf{V}}}
\nc{\sW}{{\mathsf{W}}}
\nc{\sZ}{{\mathsf{Z}}}
\nc{\sfp}{{\mathsf{p}}}
\nc{\sr}{{\mathsf{r}}}
\nc{\sfb}{{\mathsf{b}}}
\nc{\sfc}{{\mathsf{c}}}
\nc{\sd}{{\mathsf{d}}}
\nc{\sg}{{\mathsf{g}}}
\nc{\sfl}{{\mathsf{l}}}

\nc{\BK}{{\bar{K}}}

\nc{\tA}{{\widetilde{\mathbf{A}}}}
\nc{\tB}{{\widetilde{\mathcal{B}}}}
\nc{\tg}{{\widetilde{\mathfrak{g}}}}
\nc{\tG}{{\widetilde{G}}}
\nc{\TM}{{\widetilde{\mathbb{M}}}{}}
\nc{\tO}{{\widetilde{\mathsf{O}}}{}}
\nc{\tU}{{\widetilde{\mathfrak{U}}}{}}
\nc{\TZ}{{\tilde{Z}}}
\nc{\tZ}{\widetilde{Z}{}}
\nc{\tx}{{\tilde{x}}}
\nc{\tbv}{{\tilde{\bv}}}
\nc{\tfP}{{\widetilde{\mathfrak{P}}}{}}
\nc{\tz}{{\tilde{\zeta}}}
\nc{\tmu}{{\tilde{\mu}}}

\nc{\td}{\ddot{\underline{d}}{}}
\nc{\tzeta}{\widetilde{\zeta}{}}
\nc{\hd}{{\widehat{\underline{d}}}}
\nc{\hG}{{\widehat{G}}}
\nc{\hBP}{\widehat{\mathbb P}{}}
\nc{\hQ}{{\widehat{Q}}}
\nc{\hsM}{\widehat{\mathsf M}{}}
\nc{\hfM}{\widehat{\mathfrak M}{}}
\nc{\hCP}{\widehat{\mathcal P}{}}
\nc{\hCR}{\widehat{\mathcal R}{}}
\nc{\hCS}{{\widehat{\mathcal S}}}
\nc{\hfZ}{\widehat{\mathfrak Z}{}}

\nc{\urho}{\underline{\rho}}
\nc{\uB}{\underline{B}}
\nc{\uC}{{\underline{\mathbb{C}}}}
\nc{\ui}{\underline{i}}
\nc{\ofP}{{\overline{\mathfrak{P}}}}

\nc{\hrho}{{\hat{\rho}}}

\nc{\unl}{\underline}
\nc{\ol}{\overline}
\nc{\one}{{\mathbf{1}}}
\nc{\two}{{\mathbf{t}}}

\nc{\Tot}{{\mathop{\operatorname{\rm Tot}}}}
\nc{\Hilb}{{\mathop{\operatorname{\rm Hilb}}}}
\nc{\CHom}{{\mathop{\operatorname{{\mathcal{H}}\it om}}}}
\nc{\defi}{{\mathop{\operatorname{\rm def}}}}
\nc{\length}{{\mathop{\operatorname{\rm length}}}}

\nc{\Cliff}{{\mathsf{Cliff}}}
\nc{\Fl}{{\mathsf{Fl}}}
\nc{\Fib}{{\mathsf{Fib}}}
\nc{\Coh}{{\mathsf{Coh}}}
\nc{\FCoh}{{\mathsf{FCoh}}}

\nc{\reg}{{\text{\rm reg}}}

\nc{\cplus}{{\mathbf{C}_+}}
\nc{\cminus}{{\mathbf{C}_-}}
\nc{\cthree}{{\mathbf{C}_*}}
\nc{\Qbar}{{\bar{Q}}}

\nc{\bh}{{\bar{h}}}
\nc{\bOmega}{{\overline{\Omega}}}
\nc\tGr{\widetilde{\Gr}}

\nc{\seq}[1]{\stackrel{#1}{\sim}}
\nc\ogu{\overline{G/U}}
\nc\chlam{\check{\lam}}

\nc\St{\operatorname{St}}

\renewcommand\gg{\Gr_G}
\nc\uS{\underline{S}}
\nc\QM{\mathcal{QM}}
\begin{document}
\title{Semi-infinite Schubert varieties and quantum $K$-theory of flag manifolds}
\author{Alexander Braverman and Michael Finkelberg}

\begin{abstract}
Let $\fg$ be a  semi-simple Lie algebra  over $\CC$ and
let $\CB_{\fg}$ be its flag variety. In this paper we study the spaces $Z^{\alpha}_{\fg}$
of based quasi-maps $\PP^1\to \CB_{\grg}$ (introduced in \cite{FM}) as well as their affine
versions (corresponding to $\grg$ being untwisted affine algebra) introduced in \cite{BFG}.
The purpose of this paper is two-fold. First we study the singularities of the above spaces (as was explained
in \cite{FM} and \cite{Bicm} they are supposed to model singularities of the not rigorously defined ``semi-infinite Schubert varieties").
 We show that $Z^{\alp}_{\fg}$ is normal and when $\grg$ is simply laced, $Z^{\alp}_{\fg}$ is Gorenstein and has rational
singularities; some weaker results are proved also in the affine case.

The second purpose is to study the character of the
ring of functions on $Z^{\alp}_{\fg}$. When $\grg$ is finite-dimensional and simply laced we show that the generating
function of these characters satisfies
the ``fermionic formula" version of quantum difference Toda equation (cf. \cite{ferm}), thus extending the results for $\fg=\mathfrak{sl}(N)$ from
\cite{GL} and
\cite{BF0}; in view of the first part this also proves a conjecture from \cite{GL} describing the
quantum $K$-theory of $\CB_{\grg}$ in terms of the Langlands dual quantum group $U_q(\check\fg)$
(for non-simply laced $\fg$ certain modification of that conjecture is necessary). Similar analysis (modulo certain assumptions)
is performed for affine $\fg$, extending the results of \cite{BF0}.
\end{abstract}
\maketitle
\sec{int}{Introduction}
\ssec{}{Spaces of quasi-maps}
Let $G$ be a semi-simple simply connected group over $\CC$ with Lie algebra $\fg$; we shall
denote by $\check\fg$ the Langlands dual algebra of $\fg$. Let also $\CB_{\fg}$ denote its flag variety.
We have $H_2(\CB_{\fg},\ZZ)=\Lam$, the coroot lattice of $\fg$.
We shall denote by $\Lam_+$ the sub-semigroup
of positive elements in $\Lam$.

Let $\bfC\simeq \PP^1$ denote a (fixed) smooth connected projective curve (over $\CC$) of genus $0$; we are going to fix
a marked point $\infty\in \bfC$. For each $\alp\in\Lam_+$ we can consider the space $\calM^{\alp}_{\fg}$
of maps $\bfC\to \CB_{\fg}$ of degree $\alp$. This is a smooth quasi-projective variety.
It has a compactification $\QM^{\alp}_{\fg}$ by means of the space of {\em quasi-maps} from $\bfC$ to $\CB_{\fg}$ of degree
$\alp$. Set-theoretically this compactification can be described as follows:
\eq{strat}
\QM^{\alp}_{\fg}=\bigsqcup\limits_{0\leq \beta\leq \alp}\calM^{\beta}_{\fg}\x\Sym^{\alp-\beta}(\bfC)
\end{equation}
where $\Sym^{\alp-\beta}(\bfC)$ stands for the space of ``colored divisors" of the form
$\sum \gam_i x_i$ where $x_i\in \bfC$, $\gam_i\in \Lam_+$ and $\sum \gam_i=\alp-\beta$.

Let us fix a pair of opposite Borel subgroups $B,B_-\subset G$;
then we can write $\CB_{\fg}=G/B$. We can now
consider the space $\overset{\circ}Z{}^{\alp}_{\fg}$ of {\em based maps}
$(\bfC,\infty)\to (\CB_{\fg}=G/B,e_-)$ (here $e_-\in G/B$ denotes the
class of $B_-$, and
a map $f:\bfC\to \CB_{\fg}$ is called based if $f(\infty)=e_-$). This is a quasi-affine variety; the corresponding
space $Z^{\alp}_{\fg}$ of based quasi-maps (a.k.a. Zastava space in the terminology of \cite{FM} and \cite{FFKM}) is affine.
It possesses a stratification similar to \refe{strat} but with $\bfC$ in the right hand side of \refe{strat}
replaced by $\bfC -\infty$.

The following theorem is the first main result of this paper:
\th{main1}
\begin{enumerate}
\item
For any $\grg$ and $\alp$ the schemes $Z^{\alp}_{\fg}$ and $\QM^{\alp}_{\fg}$ are normal.
\item
Assume that $\fg$ is simply laced. Then $Z^{\alp}_{\fg}$ (and $\QM^{\alp}_{\fg}$) is Gorenstein (in particular, Cohen-Macaulay) and
has canonical (hence rational) singularities.
\end{enumerate}
\eth
\ssec{}{Connection to quantum $K$-theory of $\CB_{\fg}$}
In fact, we believe that $Z^{\alp}_\fg$ must have rational singularities for all $\fg$ (not necessarily simply laced).
Let us explain the importance of this assertion. Recall, that a scheme $Z$ has rational singularities, if for
some (equivalently, for any) resolution $\pi:\tilZ\to Z$ we have $R\pi_*(\calO_{\tilZ})=\calO_Z$.
The scheme $Z^{\alp}_\fg$ has a resolution by means of the Kontsevich moduli space $M^{\alp}_{\fg}$ of stable
maps from a nodal curve $C$ of genus $0$ to $\CB_{\fg}\x \PP^1$
which have degree $(\alp,1)$ and with some analog of the ``based" condition (cf. \refs{rational-sing} for more detail).
The space $M^{\alp}_{\fg}$ is a smooth Deligne-Mumford stack which has a natural action of $T\x\CC^*$, where
$T\subset B\subset G$ is a maximal torus (here the action of $T$ comes from the fact that it acts
on $\CB_{\fg}$ preserving $e_-$ and the action of $\CC^*$ comes from the action on $\PP^1$ preserving
$\infty$). Let $\fJ_{\alp}$ be the $T\x\CC^*$-equivariant pushforward of $\calO_{M^{\alp}_{\fg}}$ to
$\Spec(\CC)$
(i.e. the character of $[R\Gam(M^{\alp}_{\fg},\calO_{M^{\alp}_{\fg}})]$ of $R\Gam(M^{\alp}_{\fg},\calO_{M^{\alp}_{\fg}})$
 with respect to $T\x \CC^*$).  This is a rational function on $T\x \CC^*$ and we are going to
write $\fJ_{\alp}=\fJ_{\alp}(z,q)$ where $z\in T, q\in \CC^*$.

It is explained in \cite{GL} (cf. also the appendix to \cite{B-inst} for the corresponding statement in cohomology (as opposed to $K$-theory)) that every $\fJ_{\alp}$ can be thought of as some generating
function of {\em genus 0 K-theoretic Gromov-Witten invariants with gravitational descendants} of degree $\alp$.
Moreover, it is shown in \cite{LP} that $\fJ_{\alp}$'s determine all genus zero K-theoretic Gromov-Witten
invariants.\footnote{In fact, in \cite{LP} the authors work with an arbitrary smooth projective variety $X$ instead of $\calB_{\fg}$.
 In this case the definition of $\fJ_{\alp}$ is similar, however technically the push-forward must be taken with respect to certain
 {\em virtual fundamental cycle} in $K$-theory. In the case when $X$ is a
homogeneous space of a linear algebraic group this reduces to the usual
push-forward.}

Thus computing $\fJ_{\alp}$ is an important problem.
\reft{main1} implies that (for simply laced $\fg$) one can replace this equivariant pushforward
with
the character $[\calO_{Z^{\alp}_\fg}]$ of the ring of polynomial functions on $Z^{\alp}_\fg$ with respect to the action of $T\times \CC^*$.

It is often convenient to organize all $\fJ_{\alp}$
into a generating function:
 \eq{generating}
J_{\fg}(z,x,q)=\sum\limits_{\alp\in\Lam_+}x^{\alp}\fJ_{\alp},
\eeq
where $x$ lies in the dual torus
$\check T$. This function is called the {\em equivariant $K$-theoretic $J$-function} of $\calB_{\fg}$ (once again, it can be defined for any smooth projective variety $X$).
\ssec{}{Fermionic formula}
The function $J_{\fg}$ was studied in \cite{GL} for $G=\SL(N)$ and
it was shown to be an eigen-function
of the {\em quantum difference Toda integrable system} (cf. \cite{Et}, \cite{Sev});
this result was reproved in \cite{BF0} using other methods.
It was conjectured in~\cite{GL}
that the same result should hold for any $\fg$.

It is actually easy to see that verbatim this conjecture is false when $\fg$ is not simply laced. The main purpose of the
second part of the paper is to prove the above conjecture for any simply laced algebra $\fg$. More precisely, we are going to prove the following:
\th{main2}Assume that $\fg$ is simply laced. Then
the functions $\fJ_{\alp}$ satisfy the following recursive relation:
\eq{fermi-int}
\fJ_\alpha=\sum_{0\leq\beta\leq\alpha}\frac{q^{(\beta,\beta)/2}z^{\beta^*}}
{(q)_{\alpha-\beta}}\fJ_\beta.
\eeq
Here $\beta\mapsto\beta^*$ stands for the natural isomorphism between the coroot lattice of
$\fg$ and its root lattice.
\eth
The equation \refe{fermi-int} appears in \cite{ferm}, where the authors show that
\refe{fermi-int} holds precisely if and only if the generating function of the $\fJ_{\alp}$'s
is an eigen-function of the above-mentioned quantum difference Toda system. Thus
\reft{main2} and the main result of \cite{ferm} imply the following:
\cor{toda}
Let $\fg$ be simply laced. Then
the equivariant $K$-theoretic $J$-function of $\CB_{\fg}$ is an eigen-function of the quantum difference Toda
integrable system associated with $\fg$.
\ecor
One may ask whether the assumption that $\fg$ is simply laced is really essential. As was mentioned above
verbatim \reft{main2} (and thus also \refc{toda}) does not hold for non-simply laced $\fg$ (the reason for this is explained
in \refs{non-simply-laced}). On the other hand
in \refs{non-simply-laced} we show how to modify the geometric problem a little (using the recent result of
\cite{Z}) in order to make a correct statement for all $\fg$. It is worthwhile to note that the corresponding analog
of $J_{\fg}$ in that case becomes an eigen-function of the quantum difference Toda system associated with
$\check \fg$ (in the simply laced case we have $\fg=\check \fg$). The reader should compare this statement with the main
result of \cite{Kim} which deals with the ``usual" (i.e. cohomological) $J$-function of $\CB_{\fg}$.
\ssec{}{Representation-theoretic interpretation}
In this subsection we discuss possible interpretation of the above results in terms of geometric representation theory;
this subsection will not be used in the future, so uninterested reader may skip this discussion and go to \refss{norm-idea}.

\refc{toda} and the constructions of \cite{Sev} and \cite{ferm}
imply also the following:
\cor{whit}
In the simply laced case the function $J_{\fg}$ is equal to the Whittaker matrix coefficient in the universal
Verma of $U_q(\check\fg)$.
\ecor
 In \cite{BF0} this result was proved directly for $G=\SL(N)$. Namely, in that case
the space $Z^{\alp}_\fg$ has a small resolution of singularities (usually called Laumon's resolution) which we shall denote by
$\CP^{\alp}$. In \cite{BF0} we construct an action of the quantum group $U_q(\mathfrak{sl}(N))$ on
$\calV=\oplus_{\alp} K_{T\x\CC^*}(\CP^{\alp})_{\on{loc}}$
(here the subscript ``loc" means ``localized equivariant $K$-theory'') and identify the
corresponding $U_{q}(\mathfrak{sl}(N))$ module with the {\em universal Verma module}. Moreover, the natural pairing
on $\calV$ gets identified with the Shapovalov form on the Verma module. In addition if we denote
by $1_{\alp}\in K_{T\x\CC^*}(\CP^{\alp})$  then the formal sum
$\sum_{\alp} 1_{\alp}$ (lying in some completion of $\calV$) is the Whittaker vector in $\calV$ (i.e. an eigen-vector
of the positive part of $U_q(\mathfrak{sl}(N))$). It is easy to see that these results imply \refc{whit} (we refer the reader
to \cite{BF0} for the details).

It would be very interesting to prove \refc{whit} along similar lines, however we don't know how to do this, since
for general $\fg$ there is no resolution of $Z^{\alp}_\fg$ similar to $\CP^{\alp}$. In addition we would like to mention that
the notion of Whittaker vector for $\fg$ ( or $\check\fg$), which is developed in \cite{Sev} (cf. \cite{Et} for a closely related
approach) depends on a choice of orientation of the Dynkin diagram of $\fg$;
it would be very interesting to understand how it can be incorporated in the above constructions
(for $\fg=\mathfrak{sl}(N)$ there is a natural choice of orientation).


\ssec{norm-idea}{Idea of the proof of normality}Let us now go back and explain the idea of the proof of \reft{main1}(1), since
in our opinion this proof is of independent interest.

Let $\Gr_G=G((t))/G[[t]]$ be the affine Grassmannian of $G$. It is well-known that the orbits of $G[[t]]$ on $\Gr_G$
are in one-to-one correspondence with the elements of the dominant cone
$\Lam^+$; for each $\lam\in\Lam^+$ we shall denote
the corresponding orbit by $\Gr_G^{\lam}$. Its closure $\oGr_G^{\lam}$ is the union of all $\Gr_G^{\mu}$
with $\mu\leq \lam$. It is well-known (cf. e.g. \cite{Fa}) that $\oGr^{\lam}_G$ is normal, Cohen-Macaulay and has
rational singularities (in fact, it is also Gorenstein - cf. \cite{BK}).

The schemes $Z^{\alp}_\fg$ were originally defined in \cite{FM} in order to give a model for the singularities
of $\oGr_G^{\lam}$ at a point of $\Gr_G^{\mu}$ when both $\lam$ and $\mu$ are very large and $\lam-\mu=\alp$.
However, although this statement was used as a guiding principle
in many works related to $Z^{\alp}_\fg$
(cf. \cite{Bicm} for a review), it was never given any  precise meaning.

The purpose of \refs{normality} is to formulate some version of the above principle precisely. This formulation immediately
implies normality of $Z^{\alp}_\fg$ but other parts of \reft{main1} still have to be proven by other means. Roughly speaking
we show the following. Given $\lam$ and $\mu$ as above one can construct certain transversal slice $\ocalW^{\lam}_{G,\mu}$ to $\Gr_G^{\mu}$ in $\oGr_G^{\lam}$.
This transversal slice is also acted on by $T\x\CC^*$.
In \refs{normality} we construct a $T\x\CC^*$-equivariant map $\ocalW^{\lam}_{G,\lam-\alp}\to Z^{\alp}_\fg$ and we show that this map
induces an isomorphism on functions of given homogeneity degree with respect to $\CC^*$ when $\lam$ is very large.
This easily implies that $Z^{\alp}_\fg$ is normal.
\ssec{}{Affine case}The definition of the schemes $Z^{\alp}_\fg$ was generalized in \cite{BFG} to the case when $\fg$ is an untwisted
affine algebra. We conjecture that \reft{main1} and \reft{main2} hold  in this case; this should be useful for
studying the Nekrasov  partition function of 5-dimensional pure gauge theory compactified on $S^1$ (cf. \cite{Nek})
 in the spirit of \cite{BrEt}. In \refs{normal-sln} we prove
\reft{main1} for $\fg=\mathfrak{sl}(N)_{\aff}$; this easily implies  \reft{main2} in this case, in view of the results of
\cite{BF0}.
\ssec{}{Contents}This paper is organized as follows.
In \refs{normality} we discuss the relation between
$Z^{\alp}_\fg$ and the transversal slices $\ocalW^{\lam}_{G,\mu}$ in the affine
Grassmannian and prove that the schemes $Z^{\alp}_\fg$ are normal.
In \refs{normal-sln} we use a different method to show that the affine
analogs of $Z^{\alp}_\fg$
are normal, Gorenstein and have rational singularities for $G=\SL(N)$.
In \refs{determinant} we study the equation of the boundary of $Z^{\alp}_\fg$;
we use it in \refs{rational-sing}
in order to prove the second part of \reft{main1}. 
\reft{main2} is proved in~\refs{fermionic}.
Finally in \refs{non-simply-laced} we explain how to extend \reft{main1}
and \reft{main2} to non-simply laced case
using the twisted affine Grassmannian studied in \cite{Z}.
\ssec{}{Acknowledgements}The idea of proving \refc{toda} using \cite{ferm} came to us during a talk by E.~Feigin (explaining the results of \cite{ferm}) in Oberwolfach in March 2009. Also, the ideas of \refs{normality} came to us during a conversation
with J.~Kamnitzer in Jerusalem in December 2010. In addition we would like to thank B.~Feigin, A.~Knutson,
A.~Kuznetsov and X.Zhu for illuminating discussions. We are grateful to an
anonymous referee whose remarks resulted in the correction and simplification
of the proof of the main~\reft{main2}.
Part of this work was completed when both authors were visiting CIRM (Luminy)
and IESC (Cargese). The authors would like to thank those institutions for
their hospitality. A.~B. was partially supported by NSF grants DMS-0854760 and
DMS-0901274.
M.F. was partially supported by the RFBR grants 12-01-00944, 12-01-33101,
13-01-12401/13,
the National Research University Higher School of Economics' Academic Fund
award No.12-09-0062 and
the AG Laboratory HSE, RF government grant, ag. 11.G34.31.0023.
This study comprises research findings from the ``Representation Theory
in Geometry and in Mathematical Physics" carried out within The
National Research University Higher School of Economics' Academic Fund Program
in 2012, grant No 12-05-0014.
\sec{normality}{Normality of Zastava spaces via transversal slices in the affine Grassmannian}
\ssec{qm-def}{Quasi-maps and Zastava spaces}
In this section we recall the definition of $\QM^{\alp}_{\fg}$ and $Z^{\alp}_{\fg}$,
cf. e.g.~\cite[Section~3]{FM}. Choose a Borel subgroup $B\subset G$
with unipotent radical $U$. The quotient $G/U$ is a quasi-affine variety and we denote by
$\ogu$ its affine closure. The torus $T=B/U$ acts on $G/U$ on the right and this action extends to $\ogu$.

Let $\calF_T$ be a $T$-bundle over $\bfC$. For every weight $\check\lam:T\to\GG_m$ of $T$ we may consider
the associated line bundle $\calF_T^{\check\lam}$ on $\bfC$. We say that $\calF_T$ has degree $\alp\in\Lam$ if
for every $\chlam$ as above the bundle $\calF_T^{\chlam}$ has degree $\la\chlam,\alp\ra$.
Then
the scheme $\underline\QM{}^{\alp}_{\fg}$ parametrizes the following data:

a) $T$-bundle $\calF_T$ on $\bfC$ of degree $-\alp$;

b) A $T$-equivariant map $\kap:\calF_T\to \bfC\x\ogu$ of fibre bundles over
$\bfC$ such that over the
generic point of $\bfC$ this map goes to $\bfC\x G/U$.

\noindent
More explicitly, giving $\kap$ is equivalent to specifying the following data:
for every dominant $\chlam$ the map we need to specify an embedding of locally free sheaves
$\kap^{\chlam}:\calF_T^{\chlam}\to \calO_{\bfC}\ten V(\chlam)$ (where $V(\chlam)$ denotes the
corresponding irreducible representation of $G$). The maps $\kap^{\lam}$ must satisfy certain
{\em Pl\"ucker relations}; we refer the reader to \cite{Bicm} for the details. It is easy to make the above into a modular problem,
which defines
$\underline\QM{}^{\alp}_{\fg}$ as a scheme, which is irreducible and projective of dimension $2|\alp|+\dim(G/B)$ (here
$|\alp|=\la\check\rho_\fg,\alp\ra$, where $\check\rho_{\fg}$ denotes the half-sum of the positive roots of $\fg$). It follows from~\cite{fema} and~\cite{drinf}
that the scheme $\underline\QM{}^{\alp}_{\fg}$ is in general {\em not reduced},
and we define the variety $\QM{}^{\alp}_{\fg}$ as 
$(\underline\QM{}^{\alp}_{\fg})_{\operatorname{red}}$.

Given $(\calF_T,\kap)$ as above let $\bfU\subset \bfC$ be the open subset of $\bfC$ over which $\kap$ lands in
$\bfC\x G/U$. If $x\in \bfU$ then we shall say that the quasi-map $(\calF_T,\kap)$ has no defect at $x$.
It is clear that $\kap$ defines a map $f:\bfU\to G/B$. The (open dense) subset of $\QM^{\alp}_{\fg}$
consisting of those quasi-maps for which $\bfU=\bfC$ is precisely the space $\calM^{\alp}_{\fg}$ of
maps $f:\bfC\to \CB_{\fg}=G/B$ of degree $\alp$.

Let us now fix another Borel subgroup
$B_-$ such that $B\cap B_-\simeq T$; then $B_-$ defines a point
$e_-\in \CB_{\fg}=G/B$. The scheme $Z^{\alp}_{\fg}$ is a locally closed subscheme of $\QM^{\alp}_{\fg}$ which
corresponds to those quasi-maps which have no defect at $\infty\in\bfC$ and such that $f(\infty)=e_-$ where
$f$ is as above. The scheme $Z^{\alp}_{\fg}$ is an affine, reduced and irreducible scheme of
dimension $2|\alp|$. The intersection $\overset{\circ}Z{}^{\alp}_{\fg}=Z^{\alp}_{\fg}\cap \calM^{\alp}_{\fg}$ is the
space of
{\em based maps} $f:(\bfC,\infty)\to (\CB_{\fg},e_-)$ (i.e. those maps which send $\infty$ to $e_-$).

The schemes $\QM_{\fg}^{\alp}$ and $Z^{\alp}_{\fg}$ possess the following stratification:
\eq{strat-2}
\QM^{\alp}_{\fg}=\bigsqcup\limits_{0\leq \beta\leq \alp}\calM^{\beta}_{\fg}\x\Sym^{\alp-\beta}(\bfC);\quad
Z^{\alp}_{\fg}=\bigsqcup\limits_{0\leq \beta\leq \alp}\overset{\circ}Z{}^{\beta}_{\fg}\x\Sym^{\alp-\beta}(\bfC-\infty).
\end{equation}
Here for any curve $X$ and $\gam\in\Lam_+$ we denote by $\Sym^{\gam}(X)$ the scheme parametrizing
``colored divisors" $\sum\gam_i x_i$ where $x_i\in X$, $\gam_i\in \Lam_+$ and $\sum \gam_i=\gam$.

\ssec{}{The affine Grassmannian}Let $\calK=\CC((t))$, $\calO=\CC[[t]]$. By the {\it affine
Grassmannian} of $G$ we will mean the quotient
$\gg=G(\calK)/G(\calO)$. It is known (cf. \cite{BD,MV}) that
$\gg$ is the set of $\CC$-points of an ind-scheme over
$\CC$, which we will denote by the same symbol.

Since $G$ is simply connected, its coweight (=cocharacter) lattice coincides with the
coroot lattice $\Lam=\Lam_G$. We will denote the cone of dominant coweights
by $\Lambda^+\subset\Lambda$. Let $\Lam^{\vee}$
denote the dual lattice (this is the weight lattice of $G$).
We let $2\check\rho_G$ denote the sum of the positive roots of $G$.

The group-scheme $G(\calO)$ acts on $\gg$ on the left and
its orbits can be described as follows.
One can identify the lattice $\Lam_G$ with
the quotient $T(\calK)/T(\calO)$. Fix $\lam\in\Lam_G$ and
let $t^\lambda$ denote any lift of $\lam$ to $T(\calK)$.
Let $\gg^{\lam}$ denote the $\go$-orbit of $t^\lambda$
(this is clearly independent of the choice of $t^\lambda$).
The following result is well-known:
\lem{gras-orbits}
\begin{enumerate}
\item
$$
\gg=\bigcup\limits_{\lam\in\Lam_G}\gg^{\lam}.
$$
\item
We have $\Gr_G^{\lam}=\Gr_G^{\mu}$ if and only if $\lam$ and $\mu$ belong
to the same $W$-orbit on $\Lam_G$ (here $W$ is the Weyl group of $G$). In particular,
$$
\gg=\bigsqcup\limits_{\lam\in\Lam^+_G}\gg^{\lam}.
$$
\item
For every $\lam\in\Lam^+$ the orbit
$\gg^{\lam}$ is finite-dimensional and its dimension is
equal to $\la\lam,2\check\rho_{\fg}\ra$.
\end{enumerate}
\elem
Let $\ogl$ denote the closure of $\gg^{\lam}$ in $\gg$;
this is an irreducible projective algebraic variety; one has
$\gg^{\mu}\subset \ogl$ if and only if $\lam-\mu$ is a sum of positive roots of
$\check G$.
\ssec{trans-finite}{Transversal slices}Consider the group
$G[t^{-1}]=G(\BC[t^{-1}])\subset G((t))=G(\BC((t)))$; let us denote by
$G_1=G[t^{-1}]_1$ the kernel of the natural (``evaluation at $\infty$") homomorphism
$G[t^{-1}]\to G$. For any $\lam\in\Lam$ let $\Gr_{G,\lam}=G[t^{-1}]\cdot t^{\lam}$. Then
it is easy to see that one has
$$
\gg=\bigsqcup\limits_{\lam\in\Lam^+}\Gr_{G,\lam}
$$

Let also $\calW_{G,\lam}$
denote the $G_1$-orbit of $t^{\lam}$.
For any $\lam,\mu\in\Lam^+$, $\lam\geq \mu$ set
$$
\Gr^{\lam}_{G,\mu}=\gg^{\lam}\cap \Gr_{G,\mu},\quad
\oGr^{\lam}_{G,\mu}=\oGr_G^{\lam}\cap \Gr_{G,\mu}
$$
and
$$
\calW^{\lam}_{G,\mu}=\gg^{\lam}\cap \calW_{G,\mu},\quad
\ocalW^{\lam}_{G,\mu}=\oGr_G^{\lam}\cap \calW_{G,\mu}.
$$
Note that $\ocalW^{\lam}_{G,\mu}$ contains the point $t^{\mu}$ in it.

Let $\BG_m=\BC^*$ act on $\Gr_G$ by loop rotation.
\lem{enum}
\begin{enumerate}
\item
The point $t^{\mu}$ is the only $\CC^*$-fixed point in $\ocalW^{\lam}_{G,\mu}$.
The action of $\CC^*$ on $\ocalW^{\lam}_{G,\mu}$ is ``repelling", i.e. for any
$w\in \ocalW^{\lam}_{G,\mu}$ we have
$\displaystyle{\lim\limits_{a\to\infty} a(w)=t^{\mu}.}$
\item
The orbit $G\cdot t^\mu$ is a connected component of the $\BC^*$-fixed point
set $\Gr_G^{\BC^*}$, isomorphic to a partial flag variety of $G$.
The action of $\BC^*$ on $\Gr_{G,\mu}$ is ``repelling'',
i.e. for any $w\in\Gr_{G,\mu}$ we have
$\displaystyle{\lim\limits_{a\to\infty} a(w)\in G\cdot t^\mu.}$
\item
There exists an open subset $\calU$ in $\Gr_G^{\mu}$ and an open embedding
$\calU\x \ocalW^{\lam}_{G,\mu}\hookrightarrow \oGr^{\lam}_G$ such that the diagram
$$
\begin{CD}
\calU\x \{t^{\mu}\} @>>> \Gr^{\mu}_{G}\x\{t^{\mu}\}\\
@VVV @VVV\\
\calU\x \ocalW^{\lam}_{G,\mu} @>>> \oGr^{\lam}_G
 \end{CD}
$$
is commutative.
In other words, $\ocalW^{\lam}_{G,\mu}$ is a transversal slice to $\Gr^{\mu}_G$ inside
$\oGr^{\lam}_G$. \qed
\end{enumerate}
\elem
\prf The first two statements are obvious, and the third one follows
from~\cite[Propositions~1.3.1 and~1.3.2]{kt}.
\epr

\ssec{}{Functions on $\calW_{G,\mu}$}Let $\CC[\ocalW^{\lam}_{G,\mu}]$ denote the ring of functions
on $\ocalW^{\lam}_{G,\mu}$ and let
$$
\CC[\calW_{G,\mu}]=\lim\limits_{\leftarrow}\ \CC[\ocalW^{\lam}_{G,\mu}]
$$
be the ring of functions on the ind-scheme $\calW_{G,\mu}$. The group $T\x\CC^*$ acts
on $\ocalW^{\lam}_{G,\mu}$ and $\calW_{G,\mu}$ and thus it acts on the corresponding
ring of functions.

For any linear algebraic group $H$, we are going to denote by $H_n$ the subgroup of $H[t^{-1}]$ consisting
of those maps $h(t)$ which are equal to the identity $e\in H$ modulo $t^{-n}$; in particular,
$H_0=H[t^{-1}]$. Also, let $R_n=\CC[t^{-1}]/t^{-n}$; for any scheme $X$ over $\CC$ we can consider
the scheme of maps $\Spec(R_n)\to X$ which (abusing slightly the notation) we shall denote
by $X(R_n)$. Also, given a $\CC$-point $x\in X$ we shall denote by $X(R_n)_{\on{based}}$ the closed
sub-scheme of based maps
$\Spec(R_n)\to X$ (i.e. those maps which send the unique $\CC$-point of $\Spec(R_n)$ to $x$).
In particular, if $H$ is an algebraic group over $\CC$ then
$H_1/H_n=H(R_n)_{\on{based}}$ (where the role of the point $x$ is played by the identity  $e\in H$).

Let $\St_{\mu}\subset G_1$ be the stabilizer of $t^{\mu}$ in $G_1$.
Thus $\calW_{\mu}=G_1/\St_{\mu}$.
\lem{stabilizer}
\begin{enumerate}
\item
Fix $n\in \ZZ_{>0}$ and let $\mu\in\Lam^+$ satisfy the following condition:
\eq{condition}
\la\mu,\check\alp\ra\geq n\ \text{for every positive root $\check\alp$ of $\fg$}.
\end{equation}
Then the image of $\St_\mu$ in $G_1/G_n=G(R_n)_{\on{based}}$ is equal to $U_-(R_n)_{\on{based}}$. In particular, we have
a natural map $\pi_{\mu,n}:\calW_{G,\mu}\to G(R_n)_{\on{based}}/ U_-(R_n)_{\on{based}}$.
\item
Assume that condition \refe{condition} is satisfied. Then for every $k<n$ the map
$\pi_{\mu,n}^*:\CC[G(R_n)_{\on{based}}/ U_-(R_n)_{\on{based}}]\to \CC[\calW_{G,\mu}]$ induces an isomorphism on functions
of homogeneity degree $k$ with respect to $\CC^*$.
\end{enumerate}
\elem
\prf
(1) is obvious, so let us prove (2). First, let us discuss some preliminary facts about the
 algebra $\CC[G_1]$. It is clear that
any regular function $F:G\to \CC$ defines a map of ind-schemes
$F_1:G_1\to \CC[t^{-1}]$ such that for any $g(t)\in G_1$ the
constant term of $F_1(g(t))$ is equal to
$F(e)$. Thus for any $i>0$ we can define the function $a_{F,i}$ on $G_1$ as the coefficient of
$t^{-i}$ in $F_1$. It is easy to see that the algebra $\CC[G_1]$ is topologically
generated by all the $a_{F,i}$. Since every $a_{F,i}$ has degree $i$ with respect to $\CC^*$, it follows
that any function of homogeneity degree $<n$ lies in the subalgebra generated by $a_{F,i}$ with $i<n$.
On the other hand, if $i<n$ then any $a_{F,i}$ is invariant under the (normal) subgroup $G_n$ of $G_1$.
Hence, any function on $G_1$ of homogeneity degree $<n$ is invariant under $G_n$.

 Let $f$ be a function on $\calW_{G,\mu}$ of homogeneity degree $k$ with respect to $\CC^*$.
Then we can think of $f$ as a function on $G_1$ which is invariant on the right under $\St_{\mu}$.
Then the above discussion shows that
$f$ is automatically (left and right) invariant under $G_n$. In addition,
since $f$ is invariant under $\St_{\mu}$
it follows from (1) that (under the condition \refe{condition}) the function $f$ comes from a function $\of$
on $G_1/G_n\cdot \St_{\mu}=G(R_n)_{\on{based}}/U_-(R_n)_{\on{based}}$.
\epr
Now we pass to the main technical result of this Section. Let $\alp\in\Lam_+$. Then we have
a natural map $Z^{\alp}_\fg\to G_1/U_{-,1}$. This map is defined as follows:
let $(\calF_T,\kap)$ be a quasi-map in $Z^{\alp}_\fg$.
The fiber of $\calF_T$ at $\infty$ is automatically trivialized.
This trivialization uniquely extends to $\bfC -\{ 0\}$ and thus we get a
based map $\bfC-\{ 0\}\to \overline{G/U_-}$.
Restricting this map to $n$-th infinitesimal neighbourhood of $\infty$ in $\bfC$  we get a natural
morphism $Z^{\alp}_\fg\to  G(R_n)_{\on{based}}/U_-(R_n)_{\on{based}}$.
\th{normal-surjective}
\begin{enumerate}
\item
Let $\lam,\mu\in \Lam^+$ such that $\lam\geq \mu$ and
 let $\alp=\lam-\mu$. Then there exists a natural
birational $T\x\CC^*$-equivariant morphism 
$s^{\lam}_{\mu}:\ocalW^{\lam}_{G,\mu}\to Z^{-w_0\alp}_{\fg}$,
such that for any $n$ satisfying \refe{condition}, the following diagram is commutative:
\eq{CD}
\begin{CD}
\ocalW^{\lam}_{G,\mu} @>s^{\lam}_{\mu}>> Z^{-w_0\alp}_\fg\\
@V\pi_{\mu,n}VV @VVV\\
 G(R_n)_{\on{based}}/U_-(R_n)_{\on{based}}@>\on{id}>>
G(R_n)_{\on{based}}/U_-(R_n)_{\on{based}}
 \end{CD}
\end{equation}
(here we use the right vertical map described above).
\item
Assume again that \refe{condition} is satisfied.
Then the map $(s^{\lam}_{\mu})^*:\CC[Z^{-w_0\alp}_\fg]\to \CC[\ocalW^{\lam}_{G,\mu}]$
induces an isomorphism on functions of degree $<n$.
\end{enumerate}
\eth
\prf
First of all, we claim that (1) implies (2). Indeed, since $s^{\lam}_{\mu}$ is birational, it follows that
$(s^{\lam}_{\mu})^*$ is injective. On the other hand, it is surjective by \refl{stabilizer}(2) in view of
\refe{CD}.

Hence it is enough to explain the construction of $s^{\lam}_{\mu}$.
We start with a modular description of $\ocalW^{\lam}_{\mu}$.
Recall that $\Gr_G$ is the ind-scheme parametrizing
a $G$-bundle $\calF_G$ on $\bfC$ together with a trivialization on $\bfC-\{0\}$.
Also, the isomorphism classes of $G$-bundles on $\bfC$ are in one-to-one correspondence
with $\Lam^+=\Lam/W$ (see~\cite{G} or~\cite{R}). This identification can be described as follows: it is obvious that $T$-bundles
on $\bfC$ are in one-to-one correspondence with elements of $\Lam$. On the other hand, it is well-known
that any $G$-bundle on $\bfC$ has a reduction to $T$. Thus we get a surjective map from $\Lam$ to isomorphism
classes of $G$-bundles on $\bfC$ and it is easy to see that two $T$-bundles on $\bfC$ give the isomorphic induced $G$-bundles
if and only if one is obtained from the other by means of twist by some $w\in  W$. The $G[t^{-1}]$-orbit $\Gr_{G,\mu}\subset\Gr_G$ (see~\refss{trans-finite})
parametrizes the $G$-bundles
of isomorphism type $W\mu$ equipped with a trivialization on $\bfC-\{0\}$.
According to~\refl{enum} we have a contraction $c:\ \Gr_{G,\mu}\to G\cdot t^\mu$,
and $\CW_{G,\mu}=c^{-1}(t^\mu)$. It remains to describe the contraction $c$
to the partial flag variety $G\cdot t^\mu$ in modular terms.

Recall that the Harder-Narasimhan flag of a $G$-bundle $\CF_G$ is a
canonical reduction of $\CF_G$ to a parabolic subgroup
$P\subset G$,~\cite{R1}. In case $\CF_G$ is of isomorphism type $W\mu$,
the corresponding parabolic subgroup is nothing but the stabilizer of $t^\mu$
in $G$.
Let $HN(\CF_G)$ be the Harder-Narasimhan flag of a $G$-bundle
$\CF_G\in\Gr_{G,\mu}$. Since $\CF_G$ is trivialized off $0\in\bC$, the fiber
of $HN(\CF_G)$ at $\infty\in\bC$ lies in the partial flag variety
$G\cdot t^\mu$.
So the value of $c(\CF_G)$ is just the fiber of the
Harder-Narasimhan flag of $\CF_G$ at $\infty\in\bC$. All in all,
$\ocalW^{\lam}_{\mu}\subset\Gr_G$ parametrizes $G$-bundles on $\bC$ equipped
with a trivialization off $0\in\bC$ with a pole of order $\leq\lambda$ at 0,
such that the isomorphism class of $\CF_G$ is $W\mu$, and the fiber of the
Harder-Narasimhan flag of $\CF_G$ at $\infty\in\bC$ is the base point
$t^\mu\in G\cdot t^\mu$. Note that the degree of the Harder-Narasimhan
flag of $\CF_G$ is $w_0\mu$, i.e. 
$\deg\CF_T^{\check\lambda}=\langle-w_0\mu,\check\lambda\rangle$ 
in notations of~\refss{qm-def}.

Now let us view the Harder-Narasimhan flag of $\CF_G\in\ocalW^{\lam}_{G,\mu}$
as a reduction $\CF_P$ of $\CF_G$ to a parabolic subgroup $P\subset G$
(the stabilizer of $t^\mu$, containing $B_-$).
Let $L$ be the Levi quotient of $P$, and let
$L'$ be the quotient of $L$ modulo center. Then $\on{Ind}_P^{L'}\CF_P$ is
trivial. Hence the standard reduction to the Borel $B$
in the fiber of $\CF_G$ at
$\infty\in\bC$ canonically extends to the reduction of
$\on{Ind}_P^{L'}\CF_P$. Thus any $\CF_G\in\ocalW^{\lam}_{G,\mu}$ is canonically
equipped with a reduction $\varkappa$ to $B$ with
the standard fiber $e_-$ at $\infty\in\bC$.

Finally, we are ready for the construction of $s^\lambda_\mu$. Given
$\CF_G\in\ocalW^{\lam}_{G,\mu}$ equipped with an isomorphism $\sigma:\
\CF_G\iso\CF_G^{\on{triv}}$ defined off $0\in\bC$, we transfer
the canonical reduction $\varkappa$ to $\CF_G^{\on{triv}}$ to obtain
a $B$-structure $\sigma(\varkappa)$ on $\CF_G^{\on{triv}}$ with a
singularity at $0\in\bC$. Twisting by $-w_0\lambda\cdot0$ 
(i.e. in notations of~\refss{qm-def} replacing 
$\CF_T^{\check\lambda}\to\CO_\bC\otimes V(\check\lambda)$ by
$\CF_T^{\check\lambda}(\langle w_0\lambda,\check\lambda\rangle\cdot0)
\to\CO_\bC\otimes V(\check\lambda)$)
we obtain a
(regular) generalized $B$-structure $\kappa$ on $\CF_G^{\on{triv}}$ with
an {\em a priori} defect at $0\in\bC$ (cf.~\cite[Section~11]{FM}). Clearly
$\kappa$ has no defect off $0\in\bC$, its value at $\infty\in\bC$ is $e_-$,
and the degree of $\kappa$ equals $-w_0\alpha=w_0\mu-w_0\lambda$. We define
$s^\lambda_\mu(\CF_G,\sigma)=\kappa$. 

The theorem is proved.
\epr

\rem{convo}
Recall the notations of~\cite[Definition~11.7]{FM}:
$\overline\Gr{}_G^\lambda\stackrel{\bp}{\leftarrow}\CG\CZ^{-\mu}_\lambda
\stackrel{\bq}{\to}Z^{-w_0\alpha}_G$ is the restricted convolution diagram
(as before, $\alpha=\lambda-\mu$). The locally closed subvariety
$\ocalW^{\lam}_{G,\mu}\subset\overline\Gr{}_G^\lambda$ is open dense in the
image of $\bp$, and the restriction of $\bp$ to (the preimage of)
$\ocalW^{\lam}_{G,\mu}$ is an isomorphism. The proof of the theorem shows
that $s^\lambda_\mu=\bq\bp^{-1}|_{\ocalW^{\lam}_{G,\mu}}$.
\erem

\cor{norma}
$Z^\alpha_\fg$ is normal.
\ecor

\prf
$\oGr_G^{\lam}$ is normal (see~\cite{Fa}); hence $\ocalW^{\lam}_{G,\mu}$ is
normal by~\refl{enum}(3). Suppose a function $f\in\BC(Z^\alpha_\fg)$ is a root
of a unitary polynomial $f^r+a_{r-1}f^{r-1}+\ldots+a_0=0$ with coefficients
in $\BC[Z^\alpha_\fg]$. We choose $n$ bigger than the degrees of the
(homogeneous components) of the coefficients $a_i,\ 1\leq i\leq r$. Now we
choose $\mu\in\Lambda^+$ satisfying~\refe{condition} and such that
$\lambda=\mu-w_0\alpha\in\Lambda^+$. Then all the coefficients $a_i$ lie in
$\BC[\ocalW^{\lam}_{G,\mu}]$. Hence
$f\in\BC(Z^\alpha_\fg)=\BC(\ocalW^{\lam}_{G,\mu})$ lies in
$\BC[\ocalW^{\lam}_{G,\mu}]$. Moreover, the degree of (the highest homogeneous
component of) $f$ is less than $n$. Hence $f\in\BC[Z^\alpha_\fg]$.
\epr
\sec{normal-sln}{Normality of affine Zastava for $G=\SL(N)$}
The purpose of this Section is to give another prove of normality of Zastava
for $G=\SL(N)$ which works also in the affine case.
\ssec{21}{Notations}
We denote by $I$ the set of simple coroots of the affine group $G_{\on{aff}}=
\SL(N)_{\on{aff}}$. For $\alpha\in\BN[I]$ we denote by $Z^\alpha$ the Drinfeld
Zastava space. In~\cite{Q},~\cite{BF3} we have constructed a normal scheme
$\fZ^\alpha$ together with a morphism $\eta:\ \fZ^\alpha\to Z^\alpha$ giving
a bijection at the level of $\BC$-points. In this section we prove that
$\eta$ is an isomorphism.

Recall that $\fZ^\alpha$ is defined as the categorical quotient
$\sM^\alpha/\!/G_\alpha$ where $\sM^\alpha$ is the moduli scheme of
representations of a certain chainsaw quiver with relations $Q$ of dimension
$\alpha$.
According to~\cite[2.3--2.5]{BF3}, the stacky quotient $\sM^\alpha/G_\alpha$
is the moduli stack $\on{Perv}^\alpha(\CS_N,\CalD_\infty)$ of perverse coherent
sheaves on the Deligne-Mumford stack $\CS_N$ equipped with a framing at the
divisor $\CalD_\infty\subset\CS_N$. Let us denote by $\sg_\fZ:\
\on{Perv}^\alpha(\CS_N,\CalD_\infty)\to\fZ^\alpha$ the canonical map, and let
us denote by $\sg_Z:\ \on{Perv}^\alpha(\CS_N,\CalD_\infty)\to Z^\alpha$ the
composition of $\sg_\fZ$ with $\eta$. Let us denote by $z_0^Z\in Z^\alpha$
(resp. $z_0^\fZ\in\fZ^\alpha$) the unique point fixed by the loop rotation
action of $\BG_m$. According to~\cite[5.14]{BFG}, in order to prove that
$\eta$ is an isomorphism over the base field $\BC$, it suffices to check
that the inclusion $\sg_\fZ^{-1}(z_0^\fZ)\hookrightarrow\sg_Z^{-1}(z_0^Z)$ is
an equality. We will do this mimicking the argument of~\cite[5.16--5.17]{BFG}.
\ssec{22}{Perverse sheaves framed off the origin}
We consider the following closed substack $\on{Perv}^\alpha(\CS_N,\CS_N-0)$:
for a scheme $S$, its $S$-points are $S$-families of coherent perverse sheaves
on $\CS_N$ of degree $\alpha$ equipped with a framing on $\CS_N-0$, i.e. an
isomorphism $\CF|_{\CS_N-0}\simeq\CO_{\CS_N}\oplus\CO_{\CS_N}(-\CalD_0)\oplus
\ldots\oplus\CO_{\CS_N}((1-N)\CalD_0)$, and satisfying the following condition:
for any choice of presentation $h_c:\ \CS_N=\bC\times\CX_N$ and the
corresponding factorization map $\pi_c:\ Z^\alpha\to(\bC-\infty_\bc)^\alpha$
the composition $S\to\on{Perv}^\alpha(\CS_N,\CalD_\infty)\to Z^\alpha
\to(\bC-\infty_\bc)^\alpha$ sends $S$ to the point $\alpha\cdot0_\bC$.

Note that the choices of $h_c$ are parametrized by $c\in\BA^1$ according to
the choices of vertical directions ${\mathbf d}_v$ in~\cite[5.14]{BFG}.
Similarly to~\cite[Proposition~5.16]{BFG} we have
\lem{denis}
The composition $\on{Perv}^\alpha(\CS_N,\CS_N-0)\to
\on{Perv}^\alpha(\CS_N,\CalD_\infty)\to\fZ^\alpha$ is the constant
map to the point $z_0^\fZ\in\fZ^\alpha$.
\elem
\prf
For a collection $(A_l,B_l,p_l,q_l)_{l\in\BZ/N\BZ}$ representing a point of
$\on{Perv}^\alpha(\CS_N,\CalD_\infty)$, let us denote by $T_{W_l}$ any
endomorphism of the line $W_l$ obtained by composing the maps $A_k,B_i,p_j,q_r$,
and by $T_{V_l}$ any similarly obtained endomorphism of $V_l$.
It is well known that the ring of regular functions on $\fZ^\alpha$ is
generated by all the possible $T_{W_l}$'s and the traces of all the possible
$T_{V_l}$'s.

Let $\CF$ be an $S$-point of $\on{Perv}^\alpha(\CS_N,\CS_N-0)$. For an integer
$m$, let $\CF'$ be the constant $S$-family of coherent perverse sheaves on
$\CS_N$ corresponding to the torsion-free sheaf
$\fm_0^m\oplus\fm_0^m(-\CalD_0)\oplus\ldots\oplus\fm_0^m((1-N)\CalD_0)$
where $\fm_0$ is the maximal ideal of the point $0\in\CS_N$. Then, when $m$
is large enough, we can find a map $\CF'\to\CF$ respecting the framings of
both sheaves on $\CS_N-0$. The cone of this map is set-theoretically supported
at $0\in\CS_N$ and has cohomology in degrees $0,1$.

Let $(V_l,W_l,A_l,B_l,p_l,q_l)_{l\in\BZ/N\BZ}$ (resp.
$(V'_l,W'_l,A'_l,B'_l,p'_l,q'_l)_{l\in\BZ/N\BZ}$) be the linear algebra data
corresponding to $\CF$ (resp. $\CF'$). From the constructions
of~\cite[2.3--2.5]{BF3} it follows that there are maps $V'_l\to V_l,\
W'_l\iso W_l$ which commute with all the homomorphisms. Moreover, $q'_l\equiv0$.
From this we obtain that all the $T_{W_l}$'s vanish, and the only nonzero
$T_{V_l}$'s are matrices of the form
$A_{l_1}^{k_1}\circ B^{k_2}\circ A_{l_3}^{k_3}\circ\ldots\circ B^{k_{M-1}}
\circ A_{l_M}^{k_M}$ where $l_1=l_M=l$, and $B^{k_{2i}}$ stands for the composition
of $k_{2i}$ successive (composable) matrices of the form $B_r$, and
$l_{2i+1}+k_{2i}=l_{2i-1}$. It remains to show that any such matrix is traceless.

As we already noted, for any matrix $T'_{V_l}$ defined as in the previous
paragraph but with certain $B_{r}$ replaced by $p_{r+1}q_r$, the trace vanishes
(due to the cyclic invariance of the trace of a product, being equal to the
trace of the corresponding endomorphism $T'_{V'_l}=0$). Using the relation
$A_{r+1}B_r-B_rA_r+p_{r+1}q_r=0$ repeatedly we see that $\on{Tr}T_{V_l}=
\on{Tr}(A_l^k(B_{l-1}\circ B_{l-2}\circ\ldots\circ B_l)^{k'})=
\on{Tr}((B_{l-1}\circ B_{l-2}\circ\ldots\circ B_l)^{k'}A_l^k)$ for certain
$k,k'$. Therefore, it suffices to show that the characteristic polynomial
of a matrix $A_l+cB_{l-1}\circ B_{l-2}\circ\ldots\circ B_l$ equals $t^{\dim V_l}$
for all
$c\in\BC$. However, this characteristic polynomial is nothing but the value
at our point $(A_l,B_l,p_l,q_l)_{l\in\BZ/N\BZ}\in\sM^\alpha$ of the $l$-th
component of the factorization map $\sM^\alpha\to\fZ^\alpha\to Z^\alpha
\stackrel{\pi_c}{\longrightarrow}\BA^\alpha\to\BA^{(a_l)}$ (here $a_l$ is the
$l$-the component of $\alpha\in\BN[I]$). This completes the proof of the
lemma.
\epr
\lem{dennis}
For a scheme $S$, any $S$-point of the stack $\sg_Z^{-1}(z_0^Z)$ factors through
an $S$-point of $\on{Perv}^\alpha(\CS_N,\CS_N-0)$.
\elem
\prf Repeats the argument of~\cite[5.17]{BFG}. One has only to replace the
word ``trivialization'' in {\em loc. cit.} by ``framing'', and $\bS$ by $\CS_N$.
\epr

Now an application of~\cite[Lemma~5.15]{BFG} establishes
\th{gai}
Over the base field $\BC$, the morphism $\eta:\ \fZ^\alpha\to Z^\alpha$
is an isomorphism.
\eth
\cor{gait}
Over the base field $\BC$, the Zastava scheme $Z^\alpha$ is reduced, normal,
Gorenstein, and has rational singularities.
\ecor

\prf
$\fZ^\alpha$ is proved to be reduced and normal in~\cite[Theorem~2.7]{Q}.
Recall the resolution $\varpi:\ \CP^\alpha\to Z^\alpha$ by the affine Laumon
space (see e.g. {\em loc. cit.}). We will prove that $\CP^\alpha$ is Calabi-Yau.
Then it follows by the Grauert-Riemenschneider Theorem that $Z^\alpha$ has
rational singularities and is Gorenstein. In order to prove
that $\CP^\alpha$ is Calabi-Yau, note that if the support of $\alpha$ is
not the whole of $I=\BZ/N\BZ$, then we are in the finite (as opposed to affine)
situation, and the Calabi-Yau property of Laumon resolution is proved
in~\cite[Theorem~3]{GL} (cf.~\cite[Corollary~4.3]{FR}).
If the support of $\alpha$ is full, recall the boundary divisor
$\partial Z^\alpha\subset Z^\alpha$ introduced in~\cite[11.8]{BFG}.
The proof of~\cite[Theorem~11.9]{BFG} shows that for an integer $M\in\BN$
the divisor $M\partial Z^\alpha$ (i.e. all the components
$\partial_l Z^\alpha,\ l\in I$, of the boundary enter with the same
multiplicity $M$) is a principal divisor. Let us denote
$\eta^{-1}(\partial_l Z^\alpha)$ by $\partial_l\CP^\alpha$. We see that
$M\sum_{l\in I}\partial_l\CP^\alpha$ is a principal divisor in $\CP^\alpha$.

Now recall the meromorphic symplectic form $\Omega$ on $\CP^\alpha$
(see~\cite[3.1--3.2]{Q}). Let $\omega=\Lambda^{\on{top}}\Omega$ be the
corresponding meromorphic volume form on $\CP^\alpha$. The calculation
of~\cite[Proposition~3.5]{Q} shows that the divisor of poles of $\omega$
equals $\sum_{l\in I}\partial_l\CP^\alpha$. We conclude that the
canonical class of $\CP^\alpha$ is torsion. However, the Picard group of
$\CP^\alpha$ has no torsion since $\CP^\alpha$ is cellular. This completes the
proof of the corollary.
\epr

\cor{evid}
The resolution by the affine Laumon space $\varpi:\ \CP^\alpha\to Z^\alpha$
induces an isomorphism $\varpi^*:\ \Gamma(Z^\alpha,\CO_{Z^\alpha})\iso
\Gamma(\CP^\alpha,\CO_{\CP^\alpha})$. The higher cohomology of the structure
sheaf of the affine Laumon space vanishes: $H^k(\CP^\alpha,\CO_{\CP^\alpha})=0$
for $k>0$.
\ecor


\sec{determinant}{The boundary of Zastava}
\ssec{Cartier}{The Cartier property}
In this section $\fg$ will be an arbitrary finite dimensional simple or
untwisted affine Lie algebra with coroot lattice $\Lambda_\fg$.
Let $\Lambda^+_\fg\subset\Lambda_\fg$ denote the cone of positive linear
combinations of positive simple coroots. Let $T$ be the torus with the
cocharacter lattice $\Lambda_\fg$. For $\alpha\in\Lambda^+_\fg$, the Zastava
scheme $Z^\alpha_\fg$ is constructed in~\cite[Section~9]{BFG} (under the name of
$\fU^\alpha_{G,B}$). It is a certain closure of the (smooth) scheme
$\overset{\circ}{Z}{}^\alpha_\fg$ of degree $\alpha$ based maps from
$(\BP^1,\infty)$ to the Kashiwara flag scheme $\CB_\fg$ of
$\fg$. In case $\fg=\mathfrak{sl}(N)_{\on{aff}}$ we have $Z^\alpha_\fg=Z^\alpha$
of~\refs{normal-sln}. The complement
$\partial Z^\alpha_\fg:=Z^\alpha_\fg-\overset{\circ}{Z}{}^\alpha_\fg$
(the boundary) is a quasi-effective Cartier divisor in $Z^\alpha_\fg$
according to~\cite[Theorem~11.9]{BFG}. More precisely, there is a rational
function $F_\alpha$ on $Z^\alpha_\fg$ whose lift to the normalization of
$Z^\alpha_\fg$ is regular and has the preimage of $\partial Z^\alpha_\fg$ as the
zero-divisor.

In general, the boundary $\partial Z^\alpha_\fg$ is not irreducible; its
irreducible components $\partial_{\alpha_i}Z^\alpha_\fg$ are numbered by the
simple coroots $\alpha_i$ which enter $\alpha$ with a nonzero coefficient.
The argument in~\cite[11.5--11.7]{BFG} gives the order of vanishing of
$F_\alpha$ at the generic point of $\partial_{\alpha_i}Z^\alpha_\fg$. To formulate
the answer we assume that {\em all} the simple coroots enter $\alpha$ with
nonzero coefficients (i.e. $\alpha$ has full support); otherwise the
question reduces to the similar one for a Levi subalgebra of $\fg$.

\lem{short}
If $\alpha_i$ is a short coroot, the order of vanishing of
$F_\alpha$ at the generic point of $\partial_{\alpha_i}Z^\alpha_\fg$ is 1;
if $\alpha_i$ is a long coroot, the order of vanishing of
$F_\alpha$ at the generic point of $\partial_{\alpha_i}Z^\alpha_\fg$ is the
square length ratio of a long and a short coroot (that is, 1, 2 or 3).
\elem

\prf
First let $\fg$ be simply laced. Since the restriction of the adjoint
representation to a basic ${\mathfrak{sl}}_2$-subalgebra
${\mathfrak{sl}}_2^{\alpha_i}\subset\fg$ is independent of $\alpha_i$ as an
${\mathfrak{sl}}_2$-module, the argument of~\cite[Proof of~Theorem~11.6]{BFG}
shows that the order of vanishing of $F_\alpha$ along each boundary component
$\partial_{\alpha_i}Z^\alpha_\fg$ is the same, namely 1.

If $\fg$ is not simply laced, we realize it as the folding of a simply laced
$\tilde\fg$, i.e. invariants of a pinning-preserving automorphism $\sigma:\
\tilde\fg\iso\tilde\fg$. We have $\Lambda_+(\fg)=\Lambda_+(\tilde\fg)^\sigma$,
and given $\alpha\in\Lambda_+(\fg)$ we will denote the corresponding
$\sigma$-invariant element of $\Lambda_+(\tilde\fg)^\sigma$ by $\tilde\alpha$.
Then $Z^\alpha_\fg\simeq(Z^{\tilde\alpha}_{\tilde\fg})^\sigma$, and
$F_\alpha=F_{\tilde\alpha}|_{Z^\alpha_\fg}$. Now if $\alpha_i$ is short, $\tilde\alpha_i$
is a $\sigma$-invariant simple coroot of $\tilde\fg$, say $\beta_i$, 
and generically $\partial_{\alpha_i}Z^\alpha_\fg$ is the
transversal intersection of $\partial_{\beta_i}Z^{\tilde\alpha}_{\tilde\fg}$
with $Z^\alpha_\fg$. If $\alpha_j$ is long, $\tilde\alpha_j$ is a sum of
$\frac{(\alpha_j,\alpha_j)}{(\alpha_i,\alpha_i)}(=2$ or 3) simple
coroots, say $\beta_j,\beta_{j'},\beta_{j''}$.
They are all disjoint in the Dynkin diagram of $\tilde\fg$, and the
intersection
$\partial_{\beta_j}Z^{\tilde\alpha}_{\tilde\fg}\cap
\partial_{\beta_{j'}}Z^{\tilde\alpha}_{\tilde\fg}
\cap\partial_{\beta_{j''}}Z^{\tilde\alpha}_{\tilde\fg}$
is generically transversal. Moreover, each of
$\partial_{\beta_j}Z^{\tilde\alpha}_{\tilde\fg},
\partial_{\beta_{j'}}Z^{\tilde\alpha}_{\tilde\fg},
\partial_{\beta_{j''}}Z^{\tilde\alpha}_{\tilde\fg}$
is generically transversal to
$Z^\alpha_\fg\subset Z^{\tilde\alpha}_{\tilde\fg}$, and generically
$\partial_{\alpha_i}Z^\alpha_\fg=Z^\alpha_\fg\cap
\partial_{\beta_j}Z^{\tilde\alpha}_{\tilde\fg}=Z^\alpha_\fg\cap
\partial_{\beta_{j'}}Z^{\tilde\alpha}_{\tilde\fg}=
Z^\alpha_\fg\cap\partial_{\beta_{j''}}Z^{\tilde\alpha}_{\tilde\fg}=
Z^\alpha_\fg\cap
\partial_{\beta_j}Z^{\tilde\alpha}_{\tilde\fg}\cap
\partial_{\beta_{j'}}Z^{\tilde\alpha}_{\tilde\fg}
\cap\partial_{\beta_{j''}}Z^{\tilde\alpha}_{\tilde\fg}$. The lemma follows.
\epr

\ssec{deg}{The degree of $F_\alpha$}
The function $F_\alpha$ is an eigenfunction of the torus $T\times\BG_m^*$.
Here $T$ (the Cartan torus) acts on $Z^\alpha_\fg$ via the change of
framing at infinity, and $\BG_m^*$ (loop rotations) acts on the source
$(\BP^1,\infty)$, and hence on $Z^\alpha_\fg$ via the transport of structure.
We denote the coordinates on $T\times\BG_m^*$ by $(z,q)$.
We define an isomorphism
$\alpha\mapsto\alpha^*$ from the coroot lattice of $(G,T)$ to
the root lattice of $(G,T)$ in the basis of simple coroots as follows:
$\alpha_i^*:=\check{\alpha}_i$ (the corresponding simple root).
For an element $\alpha$ of the coroot lattice of $(G,T)$
we denote by $z^{\alpha^*}$ the corresponding character of $T$.
\prop{obv} The eigencharacter of $F_\alpha$ is $q^{(\alpha,\alpha)/2}z^{\alpha^*}$.
\eprop

The proposition will be proved in~\refss{obvi}.

\ssec{Deligne}{Deligne pairing}
In order to compute the eigencharacter of $F_\alpha$ we recall the construction
of $F_\alpha$ following Faltings~\cite{Fa}. To this end recall that
given a family $f:\calX\to S$ of smooth projective curves and two line bundles
$\calL_1$ and $\calL_2$ on $\calX$ Deligne defines a line bundle
$\la \calL_1,\calL_2\ra$ on $S$. In terms of determinant bundles the definition is simply
\eq{deligne}
\la \calL_1,\calL_2\ra=\det Rf_*(\calL_1\ten\calL_2)\ten\det Rf_*(\calO_{\calX})\ten
(\det Rf_*(\calL_1)\ten \det Rf_*(\calL_2))^{-1}.
\end{equation}
Deligne shows that the resulting pairing
$\Pic(\calX)\x\Pic(\calX)\to \Pic(S)$ is symmetric (obvious) and bilinear (not obvious).
\ssec{}{Determinant bundles on $\Bun_T$}
Let $T$ be a torus and for consistency let $\Lam$ (resp. $\Lam^{\vee}$) denote
its coweight (resp. weight) lattice.
Let also $(\cdot,\cdot)$ be an even pairing on $\Lam$.

Let also $C$ be a smooth projective curve (say, over the field $\BC$)
and let $\Bun_T$ denote the moduli stack of
$T$-bundles on $C$. Then to the above data one associates a line bundle
$\calD$ on $\Bun_T$ in the following way. Let $e_1,...,e_n$ be a basis of $\Lam$ and let
$f_1,...,f_n$ be the dual basis. For every $i=1,...,n$ let $\calL_i$ denote the
line bundle on $\Bun_T\x C$ associated to the weight $f_i$.
Let also $a_{ij}=(e_i,e_j)\in\ZZ$. Then we define
\eq{detformula}
\calD=(\bigotimes\limits_{i=1}^n \la\calL_i,\calL_i\ra^{\ten \frac{a_{ii}}{2}})\otimes
(\bigotimes\limits_{1\leq i<j\leq n} \la \calL_i,\calL_j\ra^{\ten a_{ij}}).
\end{equation}
It is easy to see that $\calD$ does not depend on the choice of the basis (here, of course, we have to use
the statement that Deligne's pairing is bilinear).
\ssec{P1}{The case $C=\bC$} Let now $C=\bC$ (the projective line).
In this case let us denote by $'\Bun(T)$ the space of $T$-bundles trivialized
at $\infty$. This a scheme isomorphic to $\Lam\x \Spec(\BC)$.
We shall denote the pull-back of $\calD$ to $'\Bun_T$ also by $'\calD$.

Recall that $\GG_m$ acts on $\bC$. If $v$ stands for a coordinate on $\GG_m$,
and $c$ stands for a coordinate on $\bC=\BP^1$, we need the action $v(c):=v^2c$.
Note that for this action any line bundle on $\bC$ can be equipped with a
$\GG_m$-equivariant structure. Comparing to the action of~\refss{deg} we have
$q=v^2$. The $\GG_m$-action extends to an action on $'\Bun_T$ and $\Bun_T$.
Since the construction of
$\calD$ is completely natural, it follows that $'\calD$ is $\GG_m$-equivariant. Since every component of $'\Bun_T$ is a point,
this equivariance is given by a character of $\GG_m$ (i.e. an integer) for every $\gam\in\Lam$.
\lem{square}The above integer is equal to $(\gam,\gam)$.
\elem
\prf
The main observation is the following. Let $\calL$ be a line bundle on $\PP^1$ of degree $n$. It is isomorphic
to $\calO(n)$ and therefore it has a unique $\GG_m$-equivariant structure such that the action of $\GG_m$ on the
fiber at $\infty$ is trivial (this makes sense since $\infty$ is a fixed point of $\GG_m$). Then we claim that with respect to this
equivariant structure
$\GG_m$ acts on $\det R\Gam(\calL)$ by the character $v\mapsto v^{n(n+1)}$. Indeed if $n\geq 0$ then $H^1(\calL)=0$ and
$H^0(\calL)$ has dimension $n+1$ with weights $0,2,...,2n$ and their sum is $n(n+1)$. If $n<0$ then $H^0(\calL)=0$ and
$H^1(\calL)$ has dimension $-n-1$ with weights $-2,-4,...,-2-2(-n-2)$ and their sum is equal to
$-2(-n-1)-(-n-1)(-n-2)=-(-n)(-n-1)=-n(n+1)$.

Let now $\calL_1$ and $\calL_2$ be two line bundles on $\PP^1$ of degrees $n_1$ and $n_2$. Then the action
of $\GG_m$ on $\la\calL_1,\calL_2\ra$ with respect to the above $\GG_m$-equivariant structure
on $\calL_1$ and $\calL_2$ corresponds to the integer
\eq{detform-action}
(n_1+n_2+1)(n_1+n_2)-(n_1+1)n_1-(n_2+1)n_2=2n_1n_2.
\end{equation}

Let now $\calF$ be a $T$-bundle of degree $\gam$ and let $\calL_i$ be the line bundle associated with $f_i\in \Lam^{\vee}$.
Then the degree of $\calL_i$ is $n_i=f_i(\gam)$. Note that $\gam=\sum n_ie_i$ and hence
$$
(\gam,\gam)=\sum\limits_{i=1}^n a_{ii}n_i^2+\sum\limits_{1\leq i<j\leq n}2a_{ij}n_in_j.
$$
Therefore according to \refe{detform-action} the action of $\GG_m$ on the fiber of $'\calD$ at $\calF$ is by the character
$v\mapsto v^m$ where
$$
m=\sum\limits_{i=1}^n \frac{a_{ii}}{2}(2n_i^2)+\sum\limits_{1\leq i<j\leq n}2a_{ij}n_in_j=(\gam,\gam).
$$
\epr

\ssec{obvi}{Determinant bundle on Zastava space}
In this section $\fg$ is an arbitrary symmetrizable Kac-Moody Lie algebra,
and $T$ is the torus with the cocharacter lattice $\Lam=\Lambda_\fg$: the coroot
lattice of $\fg$. The Kashiwara flag scheme of $\fg$ is denoted by $\CB_\fg$,
and $\overset{\circ}{Z}{}^{\gam}$ denotes the space of based maps
$(\PP^1,\infty)\to \calB_\fg$ of degree $\gamma$ (see~\cite[Theorem~18.1]{BFG}).
We let $(\cdot,\cdot)$ be the minimal even $W$-invariant form on $\Lam$.
Then we have the natural
maps $f_{\gam}:\overset{\circ}{Z}{}^{\gam}\to \Bun_T^{\gam},\
'f_{\gam}:\overset{\circ}{Z}{}^{\gam}\to\ '\Bun_T^{\gam}$.
Consider the line bundle $f_{\gam}^*\calD=\ 'f_{\gam}^*(\ '\calD)$ on
$\overset{\circ}{Z}{}^{\gam}$.

This line bundle acquires two different trivializations for the following reasons:

1) Since $'\Bun_T^{\gam}$ is just one point, the line bundle $'\calD$
is trivial there and hence its pull-back is trivial as well.

2) According to Faltings~\cite[Section~7]{Fa},
there is a trivialization of the similar bundle
on the space of all maps from any smooth projective curve $C$ to $\calB_\fg$
(Faltings proves this only for finite-dimensional
$\fg$ but his proof works word by word for any $\fg$).

Moreover, the line bundle $f_{\gam}^*\calD$ is naturally $\GG_m$-equivariant
($\GG_m$ acts on everything);
both trivializations 1) and 2) are compatible with this structure
if in 1) we let $\GG_m$ act on the
trivial bundle via the character $v\mapsto v^{(\gam,\gam)}$
(this follows from \refl{square})
and in 2) we let $\GG_m$ act trivially on the trivial bundle
(this follows from the fact that
Faltings' construction is natural with respect to everything).
Thus 1) and 2) together give us an invertible function $F_\gamma$
on $\overset{\circ}{Z}{}^{\gam}$ such that
$F_\gamma(v\phi)=v^{(\gam,\gam)}F_\gamma(\phi)$, that is
$F_\gamma(q\phi)=q^{(\gam,\gam)/2}F_\gamma(\phi)$.

This completes the proof of the main part of~\refp{obv}: we have found the
eigencharacter of $F_\alpha$ with respect to the loop rotations. It remains
to compute the eigencharacter of $F_\alpha$ with respect to the Cartan torus.
We must check the following. Given a one-parametric subgroup
$\beta:\ \BG_m\to T$, and a general point $\phi\in Z^\alpha_\fg$,
the action map $\BG_m\to Z^\alpha_\fg,\ c\mapsto \beta(c)\cdot\phi$,
extends to the
map $a:\ \BG_m\subset\BA^1\to Z^\alpha_\fg$. We must check that the function
$F_\alpha\circ a$ on $\BA^1$ has the order of vanishing $\langle\beta,
\alpha^*\rangle$ at the origin. By the factorization property of $Z^\alpha_\fg$
the question reduces to the case $\fg={\mathfrak{sl}}(2),\ \alpha=1$, which
is obvious. It follows in particular that the function $F_\alpha$ vanishes
along the boundary of $Z^\alpha_\fg$. \qed



\sec{rational-sing}{Gorenstein property of Zastava}

Suppose $\fg$ is a simply laced simple Lie algebra.

\prop{Gor}
$Z^\alpha_\fg$ is a Gorenstein (hence, Cohen-Macaulay) scheme with
canonical (hence rational) singularities.
\eprop

\prf
We are going to apply Elkik's criterion~\cite{E} in order to prove that
$Z^\alpha_\fg$ has rational singularities. To this end we will use the
Kontsevich resolution $\pi:\ M^\alpha_\fg\to Z^\alpha_\fg$
(see~\cite[Section~8]{FFKM}). We will show that the discrepancy of $\pi$ is
strictly positive, that is the singularities of $Z^\alpha_\fg$ are canonical,
hence rational.

Recall that $\overline{M}_{0,0}(\BP^1\times\CB_\fg,(1,\alpha))$ is the moduli
space of stable maps from curves of genus zero without marked points of degree
$(1,\alpha)$ to $\BP^1\times\CB_\fg$. It is a smooth Deligne-Mumford stack
equipped with a birational projection to the space of Drinfeld quasimaps from
$\BP^1$ to $\CB_\fg$.
If $C$ is such a curve of genus 0, then it has a
distinguished irreducible component $C_h$ ($h$ for horizontal)
which maps isomorphically onto
$\BP^1$. Using this isomorphism to identify $C_h$ with $\BP^1$ we obtain
the points $0,\infty\in C_h$. Now
$M^\alpha_\fg\subset\overline{M}_{0,0}(\BP^1\times\CB_\fg,(1,\alpha))$ is the
locally closed substack cut out by the open condition that $\infty\in C_h$
is a smooth point of $C$, and by the closed condition that the stable map
$\phi:\ C\to\BP^1\times\CB_\fg$ takes $\infty\in C_h\subset C$ to the marked
point $e_-\in\CB_\fg$.

The open substack $\overset{\circ}{M}{}^\alpha_\fg\simeq
\overset{\circ}{Z}{}^\alpha_\fg$ of genuine based maps is formed by the pairs
$(C,\phi)$ such that $C$ is irreducible. The complement is a normal crossing
divisor with irreducible components $D_\beta$ numbered by all
$\beta\leq\alpha$. The generic point of $D_\beta$ parametrizes the pairs
$(C,\phi)$ such that $C=C_h\cup C_v$ is a union of 2 irreducible components
($v$ for vertical), and the degree of $\phi|_{C_v}$ equals $(0,\beta)$.
If $\beta=\alpha_i$ is a simple root, then $\pi:\ D_{\alpha_i}\to Z^\alpha_\fg$
is a birational isomorphism onto the boundary component
$\partial_{\alpha_i}Z^\alpha_\fg\subset Z^\alpha_\fg$. If $\beta$ is not a simple
root, then $D_\beta$ is an exceptional divisor of the Kontsevich resolution
$\pi:\ M^\alpha_\fg\to Z^\alpha_\fg$.

Recall the symplectic form $\Omega$ on $\overset{\circ}{Z}{}^\alpha_\fg$
constructed in~\cite{fkmm}. Its top exterior power $\Lambda^{|\alpha|}\Omega$
is a nonvanishing section of the canonical line bundle on
$\overset{\circ}{Z}{}^\alpha_\fg$. It is well known that the complement
$\overset{\bullet}{Z}{}^\alpha_\fg\subset Z^\alpha_\fg$ to the union of
codimension at least 2
boundary components is smooth. Let us denote its canonical line bundle by
$\overset{\bullet}{\omega}$. Then according to~\cite[Remark~3]{fkmm},
$\Lambda^{|\alpha|}\Omega$ has poles of order 1 at all the boundary divisors
$\partial_{\alpha_i}Z^\alpha_\fg\subset\overset{\bullet}{Z}{}^\alpha_\fg$.
According to~\refl{short}, the product $F_\alpha\Lambda^{|\alpha|}\Omega$ is a
regular nowhere vanishing section of $\overset{\bullet}{\omega}$, hence
$\overset{\bullet}{\omega}$ is a trivial line bundle. We conclude that
$Z^\alpha_\fg$ is $\BQ$-Gorenstein with trivial canonical class $\omega_Z$,
and the
discrepancy of the Kontsevich resolution $\pi:\ M^\alpha_\fg\to Z^\alpha_\fg$
is isomorphic to its canonical class $\omega_M$. We have
$\omega_M\otimes\pi^*\omega_Z^{-1}=\sum_{\beta\leq\alpha}m_\beta D_\beta$.
We know that for
$\beta=\alpha_i$ a simple root the multiplicity $m_{\alpha_i}$ is 0, and we will
compute the multiplicities for all the rest $\beta$, and show that they are
all strictly positive. In fact, due to the factorization property of Zastava,
it suffices to compute a single multiplicity $m_\beta$: for $\beta=\alpha$.

\lem{desperaux}
$m_\alpha=|\alpha|+\frac{(\alpha,\alpha)}{2}-2$.
\elem

\prf
The loop rotations group $\BG_m$ (cf.~\refss{deg}) acts on $M^\alpha_\fg$ via
its action on the target $\BP^1$. The fixed point set $D_\alpha^{\BG_m}$
contains all the pairs $(C,\phi)$ such that $C$ consists of 2 irreducible
components $C_h$ and $C_v$ intersecting at the point $0\in C_h$. We will
compute $m_\alpha$ via comparison of the $\BG_m$-actions in the fibers of
$\omega_M$ and the normal bundle $\CN$ to $D_\alpha$ at such a fixed point
$(C,\phi)\in D_\alpha^{\BG_m}$.

The fiber of the normal bundle $\CN_{(C,\phi)}$ equals the tensor product
of the tangent spaces at 0 to $C_h$ and $C_v$. Hence $\BG_m$ acts on
$\CN_{(C,\phi)}$ via the character $q^{-1}$ (recall that in the normalization
of~\refss{deg} and~\refss{P1} $\BG_m$ acts on the coordinate function on
$\BP^1\simeq C_h$ via the character $q$; hence it acts on the tangent space
at 0 via the character $q^{-1}$).

According to~\cite[1.3]{K}, the tangent space $T_{(C,\phi)}D_\alpha$ is
$H^1(C,\CF^\bullet(-\infty))$ where $\CF^\bullet=\CF^0\to\CF^1$ is the following
complex of sheaves in degrees~0,1: $\CF^0$ is the sheaf of vector fields on
$C$ vanishing at 0; while $\CF^1=\phi^*\CT_{\BP^1\times\CB_\fg}$ is the pullback
of the tangent sheaf of the target. Moreover, according to {\em loc. cit.},
$H^0(C,\CF^\bullet)=H^2(C,\CF^\bullet)=0$.

It is immediate to check that $H^\bullet(C,\CF^0(-\infty))=
H^0(C,\CF^0(-\infty))$ is a
3-dimensional vector space with $\BG_m$-character $1+1+1$, hence it contributes
to $\det H^1(C,\CF^\bullet(-\infty))$ the $\BG_m$-character $1$.
It is also clear that $H^\bullet(C,\phi^*\CT_{\BP^1\times\CB_\fg}(-\infty))=
H^0(C,\phi^*\CT_{\BP^1\times\CB_\fg}(-\infty))$ is a direct sum of a
$2|\alpha|$-dimensional vector space with trivial $\BG_m$-action, and a
2-dimensional vector space with $\BG_m$-character $1+q^{-1}$. Hence it
contributes the $\BG_m$-character $q^{-1}$ to
$\det H^1(C,\CF^\bullet(-\infty))$. All in all, $\BG_m$ acts on
$\det T_{(C,\phi)}M^\alpha_\fg$ with the character $q^{-2}$, and on the fiber
of the canonical bundle $\omega_{M,(C,\phi)}$ with the character $q^2$.

Now recall that the canonical class $\omega_Z$ of $Z^\alpha_\fg$ is
trivialized by the section $F_\alpha\Lambda^{|\alpha|}\Omega$ which is an
eigensection of $\BG_m$ with the character $q^{\frac{(\alpha,\alpha)}{2}+
|\alpha|}$: this follows from~\refp{obv} and~\cite[Remark~3]{fkmm}.
Hence $\BG_m$ acts on the fiber at $(C,\phi)$ of the discrepancy line bundle
$\omega_M\otimes\pi^*\omega_Z^{-1}$ with the character
$q^{2-\frac{(\alpha,\alpha)}{2}-|\alpha|}$. It coincides with the character of $\BG_m$
in the fiber $\CN_{(C,\phi)}$ raised to the power $\frac{(\alpha,\alpha)}{2}+
|\alpha|-2$. Hence $\CO(D_\alpha)$ must enter $\omega_M\otimes\pi^*\omega_Z^{-1}$
with coefficient $\frac{(\alpha,\alpha)}{2}+|\alpha|-2$. This completes the
proof of the lemma.
\epr

We return to the proof of the proposition.
Since $m_\alpha$ is positive for nonsimple $\alpha$, the singularities of
$Z^\alpha_\fg$ are canonical, hence rational, hence Cohen-Macaulay. It remains
to prove the Gorenstein property. Let us denote by $j$ the open embedding
of $\overset{\bullet}{Z}{}^\alpha_\fg$ into $Z^\alpha_\fg$. Let us denote by
$\CalD_Z$ the dualizing sheaf of $Z^\alpha_\fg$. We have to check that the
natural map $\psi:\ \CalD_Z\to j_*\overset{\bullet}{\omega}$ is an isomorphism
(the RHS is a trivial line bundle on $Z^\alpha_\fg$).
Let $\varpi: Z^\alpha_\fg\to Y$ be a finite map to a smooth affine scheme.
Then $\CalD_Z$ is locally free over $Y$, and $\varpi_*\psi$ is an isomorphism,
hence $\psi$ is an isomorphism itself.

The proposition is proved.
\epr

\sec{fermionic}{Fermionic formulas and the boundary of Zastava}

\ssec{deep}{Structure of the boundary}
If $\alpha=\beta+\gamma$ for $\alpha,\beta,\gamma\in\Lambda^+_\fg$, then
according to~\cite[Section~10]{BFG}, we have a finite morphism
$\iota_{\beta,\gamma}:\ Z^\beta_\fg\times(\bC-\infty)^\gamma\to Z^\alpha_\fg$.
Its image (a closed reduced subscheme of $Z^\alpha_\fg$) will be denoted
by $\partial_\gamma Z^\alpha_\fg$. We will denote
$Z^\beta_\fg\times(\bC-\infty)^\gamma$
(resp. $\partial Z^\beta_\fg\times(\bC-\infty)^\gamma$) by
$\widetilde{\partial}_\gamma Z^\alpha_\fg$
(resp. $\partial\widetilde{\partial}_\gamma Z^\alpha_\fg$) for short.
The image of $\partial\widetilde{\partial}_\gamma Z^\alpha_\fg$ in
$\partial_\gamma Z^\alpha_\fg$ (a closed reduced subscheme of
$\partial_\gamma Z^\alpha_\fg$) will be denoted by
$\partial\partial_\gamma Z^\alpha_\fg$.
The union of $\partial_\gamma Z^\alpha_\fg$
over all $\gamma\in\Lambda^+_\fg$ such that $\alpha-\gamma\in\Lambda^+_\fg$
and $|\gamma|=n$ (a closed reduced equidimensional subscheme of $Z^\alpha_\fg$
of codimension $n$) will be denoted $\partial_n Z^\alpha_\fg$. Here for
$\gamma=\sum c_i\alpha_i$ we set $|\gamma|=\sum c_i$. The disjoint union of
$\widetilde{\partial}_\gamma Z^\alpha_\fg$
(resp. $\partial\widetilde{\partial}_\gamma Z^\alpha_\fg$)
over all $\gamma\in\Lambda^+_\fg$ such
that $\alpha-\gamma\in\Lambda^+_\fg$ and $|\gamma|=n$ will be denoted by
$\widetilde{\partial}_n Z^\alpha_\fg$
(resp. $\partial\widetilde{\partial}_n Z^\alpha_\fg$).
Thus we have a finite morphism $\iota_n:\ \widetilde{\partial}_n Z^\alpha_\fg\to
\partial_n Z^\alpha_\fg$, and the reduced preimage of
$\partial_{n+1}Z^\alpha_\fg\subset\partial_n Z^\alpha_\fg$ is
$\partial\widetilde{\partial}_n Z^\alpha_\fg\subset
\widetilde{\partial}_n Z^\alpha_\fg$.
Hence we have an embedding $\iota_n^*:\ \Gamma(\partial_n Z^\alpha_\fg,
\CO_{\partial_n Z^\alpha_\fg}(-\partial_{n+1}Z^\alpha_\fg))\hookrightarrow
\Gamma(\widetilde{\partial}_n Z^\alpha_\fg,
\CO_{\widetilde{\partial}_n Z^\alpha_\fg}(-\partial\widetilde{\partial}_n
Z^\alpha_\fg))$.

\ssec{K}{Equivariant $K$-theory of affine Laumon spaces}
This subsection deals with the case $\fg={\mathfrak{sl}}(N)_{\on{aff}}$.
We recall some facts from~\cite{BF0} and~\cite{ferm}. We consider the
equivariant $K$-theory of the affine Laumon spaces $\CP^\alpha$ with respect
to certain torus $\widehat{T}=\widetilde{T}\times\BC^*\times\BC^*$. Here
$\widetilde{T}$ is a certain $2^{N-1}$-fold cover of a Cartan torus
$T\subset\SL(N)$ acting on $\CP^\alpha$ via the change of framing at infinity,
while $\BC^*\times\BC^*$ acts on $\CS_N$ by dilations, and hence on $\CP^\alpha$
by the transport of structure. The coordinates on $\widehat{T}$ are denoted
by $(t_1,\ldots,t_N,u,v),\ t_1\cdot\ldots\cdot t_N=1$. Certain natural
correspondences between the affine Laumon spaces give rise to the action of
the affine quantum group $\CU$ of type $\widetilde{A}_{N-1}$
(a $\BC[v^{\pm1}]$-algebra) on $\CM=\bigoplus_\alpha\CM_\alpha:=
\bigoplus_\alpha K^{\widehat{T}}(\CP^\alpha)
\otimes_{K^{\widehat T}(pt)}\operatorname{Frac}(K^{\widehat T}(pt))$.
The Cartan subalgebra of $\CU$ is
$\CU_0=\BC[v^{\pm1},C^{\pm1},L_i^{\pm1},\ 0\leq i\leq N-1]$ in notations
of~\cite[3.5]{BF0}.
Now $\CM$ contains the universal Verma module $\CM'$ over
$\CU'=\CU\otimes\on{Frac}(\CU_0)$; here $C$ acts as $uv^N$, while for
$\alpha=\sum_{i=0}^{N-1}d_i\alpha_i$ the action of $L_i$ on $\CM_\alpha$
is given by $t_1^{-1}\cdots t_{i}^{-1}v^{d_i+i(n-i)/2}$, see~\cite[3.4,3.17]{BF0}.
A certain completion $\widehat\CM'$ contains the Whittaker vector
$\fu=\sum_\alpha\fu_\alpha$, and the dual Whittaker vector
$\fn'=\sum_\alpha\fn'_\alpha$. According to~\cite[Corollary 3.21]{BF0},
the Shapovalov scalar product
$(\fn'_\alpha,\fu_\alpha)\in\BC(t_1,\ldots,t_N,u,v)$
equals the class $\fJ_\alpha:=[R\Gamma(\CP^\alpha,\CO_{\CP^\alpha})]
\in\operatorname{Frac}(K^{\widehat T}(pt))$ up to a monomial in $t,u,v$.
Comparing with~\cite[Theorem~3.1]{ferm} we see that the collection of
rational functions $\fJ_\alpha$ is uniquely characterized by the condition
$\fJ_0=1$, and the recursion relation
\eq{fermi}
\fJ_\alpha=\sum_{0\leq\beta\leq\alpha}\frac{q^{(\beta,\beta)/2}z^{\beta^*}}
{(q)_{\alpha-\beta}}\fJ_\beta
\eeq
where $q=v^2$, and $(q)_\gamma:=\prod_{i=0}^{N-1}\prod_{s=1}^{c_i}(1-q^s)$ for
$\gamma=\sum_{i=0}^{N-1}c_i\alpha_i$;
while $z^{\gamma^*}:=\prod_{i=0}^{N-1}z_i^{c_i}$,
and $z_i=t_{i+1}t_i^{-1}u^{\delta_{0,i}}$ corresponds to the highest weight of
the standard Cartan generator
$K_i=L_i^2L_{i+1}^{-1}L_{i-1}^{-1}C^{\delta_{i,0}}\in\CU_0$
($i$ is understood as a residue mod $N$ in the latter formula).
Now~\refc{evid} yields the following

\cor{evi}
The classes $\fJ_\alpha=[\Gamma(Z^\alpha,\CO_{Z^\alpha})]\in
\operatorname{Frac}(K^{\widehat T}(pt))$ satisfy the recursion
relation~\refe{fermi}.
\ecor

\prop{clef}
Let $\fg$ be a simply laced finite or affine Lie algebra such that
$Z^\alpha_\fg$ is normal for every $\alpha$. The embedding
$\iota_n^*:\ \Gamma(\partial_n Z^\alpha_\fg,
\CO_{\partial_n Z^\alpha_\fg}(-\partial\partial_n Z^\alpha_\fg))
\hookrightarrow\Gamma(\widetilde{\partial}_n Z^\alpha_\fg,
\CO_{\widetilde{\partial}_n Z^\alpha_\fg}
(-\partial\widetilde{\partial}_n Z^\alpha_\fg))$ is an isomorphism for
any $n$ (equivalently,
$\iota_{|\gamma|}^*:\ \Gamma(\partial_\gamma Z^\alpha_\fg,
\CO_{\partial_\gamma Z^\alpha_\fg}(-\partial\partial_\gamma Z^\alpha_\fg))
\hookrightarrow\Gamma(\widetilde{\partial}_\gamma Z^\alpha_\fg,
\CO_{\widetilde{\partial}_\gamma Z^\alpha_\fg}
(-\partial\widetilde{\partial}_\gamma Z^\alpha_\fg))$ is an isomorphism
for any $\gamma\leq\alpha$),
if and only if the fermionic recursion~\refe{fermi-int} holds
for any $\alpha$. \qed
\eprop

\prf
We have $[\Gamma(\widetilde{\partial}_\gamma Z^\alpha_\fg,
\CO_{\widetilde{\partial}_\gamma Z^\alpha_\fg})]=\fJ^\fg_{\alpha-\gamma}\cdot
\frac{1}{(q)_\gamma}$ where $\fJ^\fg_{\beta}$ stands for the class
of $[\Gamma(Z^\beta_\fg,\CO_{Z^\beta_\fg})]$ in $\on{Frac}(K^{T\times\BC^*}(pt))$.
Also, $[\Gamma(\partial Z^\beta_\fg,\CO_{\partial Z^\beta_\fg})]=
(1-q^{(\beta,\beta)/2}z^{\beta^*})\fJ^\fg_\beta$ since the (reduced) subscheme
$\partial Z^\beta_\fg\subset Z^\beta_\fg$ is cut out by the equation $F_\beta$
whose $T\times\BG_m$-degree is given by~\refp{obv}. In effect, the
zero-subscheme of $F_\beta$ is generically reduced (at each irreducible
component) by~\refl{short}, and hence reduced due to normality
of $Z^\beta_\fg$. In effect, we must check that any function
$f\in\Gamma(Z^\beta_\fg,\CO_{Z^\beta_\fg})$ vanishing at the boundary is divisible
by $F_\beta$. The rational function $f/F_\beta$ is regular at the generic points
of all the components of the boundary due to~\refl{short}, so it is regular
due to normality of $Z^\beta_\fg$.

All in all we see that
$[\Gamma(\widetilde{\partial}_\gamma Z^\alpha_\fg,
\CO_{\widetilde{\partial}_\gamma Z^\alpha_\fg}
(-\partial\widetilde{\partial}_\gamma Z^\alpha_\fg))]$
is equal to
$$
[\Gamma(\widetilde{\partial}_\gamma Z^\alpha_\fg,
\CO_{\widetilde{\partial}_\gamma Z^\alpha_\fg})]-
[\Gamma(\partial\widetilde{\partial}_\gamma Z^\alpha_\fg,
\CO_{\partial\widetilde{\partial}_\gamma Z^\alpha_\fg})]=
\frac{q^{(\alpha-\gamma,\alpha-\gamma)/2}z^{(\alpha-\gamma)^*}}{(q)_\gamma}~
\fJ_{\alpha-\gamma}.
$$

Furthermore, since $[\Gamma(\partial_n Z^\alpha_\fg,
\CO_{\partial_n Z^\alpha_\fg}(-\partial_{n+1}Z^\alpha_\fg))]=
[\Gamma(\partial_n Z^\alpha_\fg,\CO_{\partial_n Z^\alpha_\fg})]-
[\Gamma(\partial_{n+1}Z^\alpha_\fg,\CO_{\partial_{n+1}Z^\alpha_\fg})]$ we have
$[\Gamma(Z^\alpha_\fg,\CO_{Z^\alpha_\fg})]=\sum_{n\geq0}
[\Gamma(\partial_n Z^\alpha_\fg,
\CO_{\partial_n Z^\alpha_\fg}(-\partial_{n+1}Z^\alpha_\fg))]=
\sum_{\gamma\leq\alpha}[\Gamma(\partial_\gamma Z^\alpha_\fg,
\CO_{\partial_\gamma Z^\alpha_\fg}(-\partial\partial_\gamma Z^\alpha_\fg))]$.
Let us view this equality as an equality of formal power series in $q,z$
with {\em nonnegative} powers and with {\em nonnegative} integral coefficients.
Note that $[\Gamma(\widetilde{\partial}_\gamma Z^\alpha_\fg,
\CO_{\widetilde{\partial}_\gamma Z^\alpha_\fg}
(-\partial\widetilde{\partial}_\gamma Z^\alpha_\fg))]\succeq
[\Gamma(\partial_\gamma Z^\alpha_\fg,
\CO_{\partial_\gamma Z^\alpha_\fg}(-\partial\partial_\gamma Z^\alpha_\fg))]$
(meaning that the LHS series is termwise bigger than or equal to the RHS series)
since $\Gamma(\partial_\gamma Z^\alpha_\fg,
\CO_{\partial_\gamma Z^\alpha_\fg}(-\partial\partial_\gamma Z^\alpha_\fg))
\hookrightarrow\Gamma(\widetilde{\partial}_\gamma Z^\alpha_\fg,
\CO_{\widetilde{\partial}_\gamma Z^\alpha_\fg}
(-\partial\widetilde{\partial}_\gamma Z^\alpha_\fg))$.

Comparing to the equality~\refe{fermi-int}, in view of the equality
$[\Gamma(\widetilde{\partial}_\gamma Z^\alpha_\fg,
\CO_{\widetilde{\partial}_\gamma Z^\alpha_\fg}
(-\partial\widetilde{\partial}_\gamma Z^\alpha_\fg))]=
\frac{q^{(\alpha-\gamma,\alpha-\gamma)/2}z^{(\alpha-\gamma)^*}}{(q)_\gamma}
\fJ_{\alpha-\gamma}$, we must
have $[\Gamma(\widetilde{\partial}_\gamma Z^\alpha_\fg,
\CO_{\widetilde{\partial}_\gamma Z^\alpha_\fg}
(-\partial\widetilde{\partial}_\gamma Z^\alpha_\fg))]=
[\Gamma(\partial_\gamma Z^\alpha_\fg,
\CO_{\partial_\gamma Z^\alpha_\fg}(-\partial\partial_\gamma Z^\alpha_\fg))]$ 
which completes the proof of the proposition.
\epr



\ssec{prth}{Proof of~\reft{main2}}




\lem{flatness}
The factorization morphism $\pi:\ Z^\alpha_\fg\to\BA^\alpha=(\bC-\infty)^\alpha$ 
is flat.
\elem

\prf
According to~\refp{Gor}, $Z^\alpha_\fg$ is Cohen-Macaulay. Evidently, 
$\BA^\alpha$ is regular. It is well-known that all the fibers of $\pi$ have
the same dimension $|\alpha|$. It remains to apply~\cite[Theorem~23.1]{M}.
\epr

Let $\CF^\alpha_\fg=\pi^{-1}(\alpha\cdot0)$ stand for the scheme-theoretic
fiber of $\pi:\ Z^\alpha_\fg\to\BA^\alpha$ over $\alpha\cdot0\in\BA^\alpha$.
Let $\fI_\alpha=[\Gamma(\CF^\alpha_\fg,\CO_{\CF^\alpha_\fg})]
\in\on{Frac}(K^{T\times\BC^*}(pt))$ be the character of the ring of regular
functions on the central fiber of $\pi$.

\cor{central fiber}
$\fI_\alpha=(q)_\alpha\fJ_\alpha$.
\ecor

\prf
First, $[\Gamma(\BA^\alpha,\CO_{\BA^\alpha})]=(q)_\alpha^{-1}$. Second, the flatness
of $\pi$ implies that $\pi_*\CO_{Z^\alpha_\fg}$ is a direct sum of finite
dimensional $T\times\BC^*$-equivariant vector bundles $\CV_\xi$ s.t. 
$T\times\BC^*$ acts in the fiber of $\CV_\xi$ over $\alpha\cdot0$ via a character
$\xi$ of $T\times\BC^*$. Finally, the fiber of $\pi_*\CO_{Z^\alpha_\fg}$ over
$\alpha\cdot0$ is nothing but $\Gamma(\CF^\alpha_\fg,\CO_{\CF^\alpha_\fg})$.
\epr

Now the fermionic recursion~\refe{fermi-int} is equivalent to
\eq{fermi-int'}
\fI_\alpha=
\sum_{\beta\leq\alpha}q^{(\beta,\beta)/2}z^{\beta^*}{\binom{\alpha}{\beta}}_q\fI_\beta
\eeq
where for $\alpha=\sum_{i\in I}a_i\alpha_i,\ \beta=\sum_{i\in I}b_i\alpha_i$ we set
${\binom{\alpha}{\beta}}_q:=
\prod_{i\in I}\prod_{s=1}^{a_i}(1-q^s)\prod_{s=1}^{b_i}(1-q^s)^{-1}\prod_{s=1}^{a_i-b_i}
(1-q^s)^{-1}$. We will prove the fermionic
recursion~\refe{fermi-int} in the equivalent formulation~\refe{fermi-int'}.
To this end, we introduce the schemes (cf. notations of~\refss{deep}): 
$\partial_\gamma\CF^\alpha_\fg,\ \widetilde\partial_\gamma\CF^\alpha_\fg,\
\partial\widetilde\partial_\gamma\CF^\alpha_\fg,\ 
\partial\partial_\gamma\CF^\alpha_\fg,\
\widetilde\partial_n\CF^\alpha_\fg,\
\partial\widetilde\partial_n\CF^\alpha_\fg$ as the scheme-theoretic fibers over
$\alpha\cdot0\in\BA^\alpha$ of the corresponding morphisms
$\partial_\gamma Z^\alpha_\fg\to\BA^\alpha$ etc. (``the fiber of $Z$ is $\CF$'').

As in the proof of~\refp{clef}, we see that
$[\Gamma(\widetilde{\partial}_\gamma\CF^\alpha_\fg,
\CO_{\widetilde{\partial}_\gamma\CF^\alpha_\fg}
(-\partial\widetilde{\partial}_\gamma\CF^\alpha_\fg))]$
is equal to
$$
[\Gamma(\widetilde{\partial}_\gamma\CF^\alpha_\fg,
\CO_{\widetilde{\partial}_\gamma\CF^\alpha_\fg})]-
[\Gamma(\partial\widetilde{\partial}_\gamma\CF^\alpha_\fg,
\CO_{\partial\widetilde{\partial}_\gamma\CF^\alpha_\fg})]=
q^{(\alpha-\gamma,\alpha-\gamma)/2}z^{(\alpha-\gamma)^*}{\binom{\alpha}{\gamma}}_q
\fI_{\alpha-\gamma}.
$$
In effect, $\widetilde\partial_\gamma\CF^\alpha_\fg$ projects to the fiber
of $\BA^{\alpha-\gamma}\times\BA^\gamma\to\BA^\alpha$ over $\alpha\cdot0$.
The character of the ring of regular functions on this fiber equals
${\binom{\alpha}{\gamma}}_q$, the projection is flat, and its fiber over
$(\alpha-\gamma)\cdot0\times\gamma\cdot0$ is nothing but $\CF^{\alpha-\gamma}_\fg$.

Furthermore, since $[\Gamma(\partial_n \CF^\alpha_\fg,
\CO_{\partial_n \CF^\alpha_\fg}(-\partial_{n+1}\CF^\alpha_\fg))]=
[\Gamma(\partial_n \CF^\alpha_\fg,\CO_{\partial_n \CF^\alpha_\fg})]-
[\Gamma(\partial_{n+1}\CF^\alpha_\fg,\CO_{\partial_{n+1}\CF^\alpha_\fg})]$ we have
$[\Gamma(\CF^\alpha_\fg,\CO_{\CF^\alpha_\fg})]=\sum_{n\geq0}
[\Gamma(\partial_n \CF^\alpha_\fg,
\CO_{\partial_n \CF^\alpha_\fg}(-\partial_{n+1}\CF^\alpha_\fg))]=
\sum_{\gamma\leq\alpha}[\Gamma(\partial_\gamma \CF^\alpha_\fg,
\CO_{\partial_\gamma \CF^\alpha_\fg}(-\partial\partial_\gamma \CF^\alpha_\fg))]$.
Let us view this equality as an equality of formal power series in $q,z$
with {\em nonnegative} powers and with {\em nonnegative} integral coefficients.
Note that $[\Gamma(\widetilde{\partial}_\gamma \CF^\alpha_\fg,
\CO_{\widetilde{\partial}_\gamma \CF^\alpha_\fg}
(-\partial\widetilde{\partial}_\gamma \CF^\alpha_\fg))]\succeq
[\Gamma(\partial_\gamma \CF^\alpha_\fg,
\CO_{\partial_\gamma \CF^\alpha_\fg}(-\partial\partial_\gamma \CF^\alpha_\fg))]$
(meaning that the LHS series is termwise bigger than or equal to the RHS series)
since $\Gamma(\partial_\gamma \CF^\alpha_\fg,
\CO_{\partial_\gamma \CF^\alpha_\fg}(-\partial\partial_\gamma \CF^\alpha_\fg))
\hookrightarrow\Gamma(\widetilde{\partial}_\gamma \CF^\alpha_\fg,
\CO_{\widetilde{\partial}_\gamma \CF^\alpha_\fg}
(-\partial\widetilde{\partial}_\gamma \CF^\alpha_\fg))$. We conclude that
the recursion~\refe{fermi-int'} holds iff
$[\Gamma(\widetilde{\partial}_\gamma \CF^\alpha_\fg,
\CO_{\widetilde{\partial}_\gamma \CF^\alpha_\fg}
(-\partial\widetilde{\partial}_\gamma \CF^\alpha_\fg))]=
[\Gamma(\partial_\gamma \CF^\alpha_\fg,
\CO_{\partial_\gamma \CF^\alpha_\fg}(-\partial\partial_\gamma \CF^\alpha_\fg))]$.

Note that both sides of the above equality are well-defined at $q=1$, 
and since we already know the inequality LHS$\succeq$RHS, it suffices to 
check the equality at $q=1$. Let us denote $\fI_\alpha|_{q=1}$ by $\fI'_\alpha$.
Then $\fI'_\alpha$ is the class of 
$[\Gamma(\ '\CF^\alpha_\fg,\CO_{'\CF^\alpha_\fg})]\in\on{Frac}(K^T(pt))$ where
$'\CF^\alpha_\fg$ is an {\em arbitrary} fiber of $\pi:\ Z^\alpha_\fg\to\BA^\alpha$.
If $\underline{x}\in\BA^\alpha$ is a general point (a configuration of 
distinct points), then by factorization $\pi^{-1}(\underline{x})$ is isomorphic 
to $\prod_{i\in I}\BA^{a_i}$, and $\fI'_\alpha=\prod_{i\in I}(1-z_i)^{-a_i}$.
The fermionic recursion~\refe{fermi-int'} at $q=1$ becomes the evident
equality $\fI'_\alpha=\sum_{\beta\leq\alpha}z^{\beta^*}\binom{\alpha}{\beta}\fI'_\beta$.
This completes the proof of~\reft{main2}.


\sec{non-simply-laced}{Non simply laced case}

\ssec{non-simp-general}{}As was mentioned before, \reft{main2} does not hold verbatim for non-simply laced $\grg$.
It is reasonable to ask whether one can modify the statement so that an
analog of \reft{main2} becomes true in the non-simply laced case.
For this one should either modify the geometric problem (i.e. modify the definition
of $\fJ_{\alp}$) or modify the equations that we want our $J$-function to satisfy.
Apparently, both ways are possible, however the former is much easier than the latter.
In this last Section we would like to explain how to modify the definition of $\fJ_{\alp}$ and sketch a proof
of the corresponding analog of \reft{main2}.
Details (as well as other variants of \reft{main2} in the non-simply laced case)
will appear in another publication.

\ssec{recursion}{Fermionic recursion}
We recall the results of~\cite{ferm}. Let $\check\fg$ be a simple Lie algebra
with the corresponding adjoint Lie group $\check G$. Let $\check T$ be
a Cartan torus of $\check G$. We choose a Borel subgroup $\check{B}\supset
\check T$. It defines the set of simple roots $\{\alpha_i,\ i\in I\}$.
Let $G\supset T$ be the Langlands dual groups. We define an isomorphism
$\alpha\mapsto\alpha^*$ from the root lattice of $(\check{G},\check{T})$ to
the root lattice of $(G,T)$ in the basis of simple roots as follows:
$\alpha_i^*:=\check{\alpha}_i$ (the corresponding simple coroot).
For two elements $\alpha,\beta$ of the root lattice of $(\check{G},\check{T})$
we say $\beta\leq\alpha$
if $\alpha-\beta$ is a nonnegative linear combination of
$\{\alpha_i,\ i\in I\}$.
For such $\alpha$ we denote by $z^{\alpha^*}$ the
corresponding character of $T$. As usually, $q$ stands for the
identity character of $\BG_m$. We set $d_i=\frac{(\alpha_i,\alpha_i)}{2}$,
and $q_i=q^{d_i}$. For $\gamma=\sum_{i\in I}c_i\alpha_i$, we set
$(q)_\gamma:=\prod_{i\in I}\prod_{s=1}^{c_i}(1-q_i^s)$.
According to~\cite[Theorem~3.1]{ferm}, the recurrence relations
$$
\calJ_\alpha=\sum_{0\leq\beta\leq\alpha}
\frac{q^{(\beta,\beta)/2}z^{\beta^*}}{(q)_{\alpha-\beta}}\calJ_\beta
$$
uniquely define a collection of rational functions $\calJ_\alpha,\ \alpha\geq0$,
on
$T\times\BG_m$, provided $\calJ_0=1$. Moreover, these functions are
nothing but
the Shapovalov scalar products of the weight components of the Whittaker
vectors in the universal Verma module over the corresponding quantum group.

We can now explain why \reft{main2} doesn't literally hold for non-simply laced $\fg$.
This can be seen, for example, in the following way. Let $\alp_i$ be a simple coroot of $\fg$. Then it is easy to see
that we have
$$
\calJ_{{\alp_i}}(z,q)= \frac{1}{(1-z^{\alp_i^*}q_i)(1-q_i)}.
$$
However, the character of $\CC[Z^{\alpha_i}_\fg]$ is simply  $\frac{1}{(1-z^{{\alp_i}^*}q)(1-q)}$.
Moreover, in this case $Z^{\alpha_i}_\fg$ is isomorphic to the affine plane $\AA^2$, so it obviously has rational singularities.
Therefore the character of $\CC[Z^{\alpha_i}_\fg]$ is equal to $\fJ_{\alp_i}$. Hence we see that if $q_i\neq q$ then
$\fJ_{\alp_i}\neq \calJ_{\alp_i}$.
\ssec{interpretation}{Geometric interpretation of $\calJ_{\alp}$}
We are now going to introduce a scheme
$\widehat{Z}{}^\alpha_\fg$ equipped with the action of $T\times\BG_m$ such
that the character of $\BC[\widehat{Z}{}^\alpha_\fg]$ equals $\calJ_\alpha$.

To this end we realize $\check\fg$ as a {\em folding} of a simple
simply laced Lie algebra
$\check\fg'$, i.e. as invariants of an outer automorphism $\sigma$ of
$\check\fg'$ preserving a Cartan subalgebra $\check\ft'\subset\check\fg'$
and acting on the root system of $(\check\fg',\check\ft')$.
In particular,
$\sigma$ gives rise to the same named automorphism of the Langlands dual
Lie algebras $\fg'\supset\ft'$. We denote by $\Xi$ the finite cyclic group
generated by $\sigma$. Let $G'\supset T'$ denote the corresponding simply
connected Lie group and its Cartan torus.
The {\em coinvariants} $X_*(T')_\sigma$ of $\sigma$ on the coroot
lattice $X_*(T')$ of $(\fg',\ft')$ coincide with the root lattice
of $\check\fg$. We have an injective map $a:\ X_*(T')_\sigma\to X_*(T')^\sigma$
from coinvariants to invariants defined as follows: given a coinvariant
$\bar\alpha$ with a representative $\alpha\in X_*(T')$ we set
$a(\bar\alpha):=\sum_{\xi\in\Xi}\xi(\alpha)$.
Given $\bar\alpha\geq0$ in the root lattice of $\check\fg$, we define an
automorphism $\varsigma$ of the based quasimaps' space $Z^{a(\bar\alpha)}_{\fg'}$
as follows. It is the composition of two automorphisms: a) $\sigma$ on the
target; b) multiplication by $\zeta$ on the source $\bC\cong\BP^1$. Here
$\zeta$ is a primitive root of unity of the order equal to the order of
$\sigma$. One can check that the fixed point set
$(\overset{\circ}{Z}{}^{a(\bar\alpha)}_{\fg'})^\varsigma$ is connected.
We define $\widehat{Z}{}^{\bar\alpha}_\fg$ as the closure of
$(\overset{\circ}{Z}{}^{a(\bar\alpha)}_{\fg'})^\varsigma$ in $Z^{a(\bar\alpha)}_{\fg'}$.

The equality $\calJ_{\bar\alpha}=[\BC[\widehat{Z}{}^{\bar\alpha}_\fg]]$ is proved
along the
lines of the argument of the previous sections. In particular, the role of
the affine Grassmannian of $G$ in the simply laced case is played by the
ramified Grassmannian of $(G',\sigma)$, see~\cite{Z}.

\bigskip
\footnotesize{
{\bf A.B.}: Department of Mathematics, Brown University,
151 Thayer St., Providence RI
02912, USA;\\
{\tt braval@math.brown.edu}}

\footnotesize{
{\bf M.F.}: IMU, IITP and National Research University
Higher School of Economics\\
Department of Mathematics, 20 Myasnitskaya st, Moscow 101000 Russia;\\
{\tt fnklberg@gmail.com}}

\end{document}